\documentclass[]{article}
\usepackage[centertags]{amsmath}
\usepackage{amssymb,amsthm,amsfonts}
\usepackage[colorlinks=true,linkcolor=blue,citecolor=blue]{hyperref}
\usepackage{url}
\usepackage{multirow}
\usepackage{enumitem}
\usepackage[normalem]{ulem}

\numberwithin{equation}{section}
\usepackage{xcolor,graphicx}

\usepackage[colorinlistoftodos,prependcaption,textsize=tiny]{todonotes}

\title{The Projected Polar Proximal Point Algorithm Converges Globally}
\author{Scott B. Lindstrom \\ Hong Kong Polytechnic University}

\def\R{\hbox{$\mathbb R$}}
\def\N{\hbox{$\mathbb N$}}

\newcommand{\argminn}[1]{\underset{#1}{\argmin}}

\renewcommand{\vec}[1]{\mathbf{#1}}
\newcommand{\PPPPA}{{$\mathbf{P^4A}$}}
\newcommand{\GPPPPA}{{$\mathbf{GP^4A}$}}
\newcommand{\Id}{{\rm Id}}
\newcommand{\cone}{{\rm cone}}

\newcommand{\dom}{{\rm dom}}
\newcommand{\zer}{{\rm zer}}
\newcommand{\lev}{{\rm lev}}
\newcommand{\Fix}{{\rm Fix}}
\newcommand{\cdk}{{\overline{\dom \kappa}}}

\newtheorem{theorem}{Theorem}[section]
\newtheorem{lemma}[theorem]{Lemma}
\newtheorem{proposition}[theorem]{Proposition}

\newtheorem{corollary}[theorem]{Corollary}
\theoremstyle{definition}
\newtheorem{definition}{Definition}

\newtheorem{example}{Example}
\newtheorem{remark}{Remark}
\newtheorem{fact}[theorem]{Fact}

\def\argmin{\mathop{\rm argmin}}

\begin{document}

\maketitle

\begin{abstract}
	Friedlander, Mac\^{e}do, and Pong recently introduced the projected polar proximal point algorithm (P4A) for solving optimization problems by using the closed perspective transforms of convex objectives. We analyse a generalization (GP4A) which replaces the closed perspective transform with a more general closed gauge. We decompose GP4A into the iterative application of two separate operators, and analyse it as a splitting method. By showing that GP4A and its under-relaxations exhibit global convergence whenever a fixed point exists, we obtain convergence guarantees for P4A by letting the gauge specify to the closed perspective transform for a convex function. We then provide easy-to-verify sufficient conditions for the existence of fixed points for the GP4A, using the Minkowski function representation of the gauge. Conveniently, the approach reveals that global minimizers of the objective function for P4A form an exposed face of the dilated fundamental set of the closed perspective transform.
\end{abstract}%Let C1={(x,y) | y > x \geq 0}. Let kappa be 0 at 0 and the indicator of C1 elsewhere. Then kappa is convex and a gauge, but it is not closed. So I must include closed.

{\small
	
}

\section{Introduction}

Friedlander, Mac\^{e}do, and Pong introduced the projected polar proximal point algorithm (\PPPPA, Definition~\ref{def:p4}) as the first proximal-point-like algorithm based on the \textit{polar envelope} of a gauge \cite{friedlander2019polar}. The motivation to study such algorithms stems from the polar envelope's relationship to infimal max convolution; Friedlander, Mac\^{e}do, and Pong showed that this relationship is analogous to the connection between the Moreau envelope and infimal convolution. They also illuminated useful variational properties of a duality framework admitted by such problems \cite{friedlander2014gauge}. It is these special properties of the algorithm, and the rich associated theoretical framework that motivated its construction, that make its global convergence an interesting question. 

The method makes use of the closed perspective transform for a proper convex function $f:X \rightarrow \R$:
\begin{equation*}
f^\pi:X\times \R_+\rightarrow \R:\;(x,\lambda)\mapsto\begin{cases}
\lambda f(\lambda^{-1}x) & \text{if}\;\;\lambda>0;\\
f_\infty(x) & \text{if}\;\;\lambda=0;\\
\infty & \text{if}\;\; \lambda < 0.
\end{cases}
\end{equation*}
Here $f_\infty$ denotes the recession function of $f$ \cite[Definition~2.5.1]{auslender2006asymptotic}, which satisfies
$$
{\rm epi}(f_\infty) = ({\rm epi}f)_\infty,
$$
and ${\rm epi}h$ denotes the epigraph of a proper convex function $h$ and $C_\infty$ denotes the recession cone of a set $C$ (see, for example, \cite[Chapter~6]{RW98}). The perspective $f^\pi$ is proper closed convex \cite[Page~67]{Roc70} and is characterized by
$$
{\rm epigraph} f^\pi = {\rm cl}\left({\rm cone}({\rm epigraph}f \times \{1\})\right).
$$
Friedlander, Mac\^{e}do, and Pong also provided a result \cite[Theorem~7.5]{friedlander2019polar} showing convergence under an assumption of strong convexity of $(f^\pi)^2$; however, the question of convergence more generally has remained open until now.

\subsubsection*{Outline and contributions}

In Section~\ref{sec:preliminaries}, we recall familiar notation and concepts from convex analysis. In Section~\ref{sec:polarprox}, we recall and analyse the polar proximity operator, a fundamental component of \PPPPA. In particular, we show that it is \textit{firmly quasinonexpansive} (Theorem~\ref{thm:TQNE}).

In Section~\ref{sec:PPPPA}, we recall the algorithm \PPPPA \; and introduce a generalization thereof, \GPPPPA  \;(Definition~\ref{def:GPPPPA}). Our motivation in so doing is that \GPPPPA \;may be described and studied as a 2-operator splitting method, a flexibility afforded by its definition on the lifted space $X \times \R$. In Section~\ref{sec:fixedpoints}, we exploit this flexibility to show that, when the operator associated to \GPPPPA \; has a nonempty fixed point set, its fixed points all share a special property (Proposition~\ref{prop:fixedpoints}), on which our analysis depends. In Section~\ref{sec:convergence}, we use this property to show that the operator is \textit{strictly quasinonexpansive} (Theorem~\ref{thm:PTQNE}) and admits global convergence of sequences to a fixed point, whenever one exists (Theorem~\ref{thm:PTstrong}). We show similar results for the algorithm's under-relaxed variants (Theorem~\ref{thm:PTrelaxedstrong}). In Section~\ref{sec:shadow}, we provide convergence results for the associated shadow sequences. In Section~\ref{sec:counterexample}, we provide an example that shows the operator associated with \GPPPPA \; is not, generically, firmly quasinonexpansive.

In Section~\ref{sec:Fixnonempty}, we provide sufficient conditions to guarantee fixed points of \GPPPPA \; (Theorem~\ref{thm:existence}). Moreover, when \GPPPPA \; specifies to \PPPPA, we show that set of global minimizers of $f$ defines an exposed face of the \textit{fundamental set} for the perspective function $f^\pi$ (Theorem~\ref{thm:fixedpoints_argmin_new}). The latter results connect to known results (\cite[Theorem~7.4]{friedlander2019polar}) about fixed points of \PPPPA, and we explain how (Remark~\ref{rem:synchronicity}).

In Section~\ref{sec:conclusion}, we connect our sufficient conditions for fixed point existence, with our convergence results that depend on that existence, to state a simple global convergence guarantee (Theorem~\ref{thm:convergence}). It shows that \PPPPA \; is globally convergent in the full generality of \cite{friedlander2019polar}.

\section{Preliminaries}\label{sec:preliminaries}

Throughout, $X$ is a finite dimensional Euclidean space with the Euclidean norm, and we will work extensively with the space $X \times \R$. For ease of clarity, when we work with a 2-tuple $(y,\lambda) \in X \times \R$, it should be understood that $y \in X$ and $\lambda \in \R$. In such a case, variable $y$ will not be bolded. In order to be succinct, we will sometimes forego the use of a 2-tuple and simply use a single variable $\vec{x} \in X \times \R$. In such a case, the bolded variable $\vec{x}$ reminds that $\vec{x} \in X \times \R$. Throughout, $\kappa:X \times \R \rightarrow \R$ is a closed \textit{gauge} in the sense of \cite{friedlander2014gauge}. In other words $\kappa$ is convex, and
$$
(\forall \vec{x} \in X\times \R)\; (\forall \lambda \in \left[0,\infty\right])\;0 \leq \kappa(\lambda \vec{x}) = \lambda \kappa (\vec{x}).
$$
For an operator $U:X\times \R \rightarrow X\times \R$, $\Fix U:=\{\vec{x} \in X \times \R \;|\; U(\vec{x})=\vec{x}\}$ is its fixed point set. For a function $f:X \rightarrow \R_+$, $\argmin f:=\{x \;|\; f(x)=\inf f(X) \}$ is its set of global minimizers, $\dom f := \{x\;|\; f(x)< \infty \}$ is its domain, $\lev_{\leq r}(f):=\{x \;|\; f(x)\leq r \}$ is its $r$-\textit{lower level set}, and $\zer f := \{x \;|\; f(x)= 0 \}$ is its zero set. For a gauge $\kappa: X \times \R \rightarrow \R_+$, the definitions of $\argmin \kappa$, $\dom \kappa$, $\lev_{\kappa < r}\kappa$, and $\zer \kappa$ are respectively analogous subsets of $X \times \R$. For a closed, convex subset $C \subset X \times \R$, $\cone C := \{\lambda \vec{x} \;|\; (\vec{x},\lambda) \in (C \times \left[0,\infty \right[) \}$ is the \textit{cone} of $C$, $P_C:\vec{x} \mapsto \argmin_{\vec{y} \in C}\|\vec{y}-\vec{x}\|$ is the projection operator associated with $C$, and $N_C(\vec{x})$ is the normal cone to $C$ at a point $\vec{x} \in C$. For projection operators associated with (closed, convex) lower level sets, we use the shorthand: $P_{f \leq r}:=P_{\lev_{\leq r}(f)}$.

We will make use of various notions of nonexpansivity, which we now introduce; more information may be found in \cite{BC}, and a comparison of what may be shown through different cutter and projection methods is found in \cite{DLR18}. The following definition may be found in either of these.

\begin{definition}[Properties of operators]\label{def:SQNE}
	Let $D\subset X\times \R$ be nonempty and let $U:D\rightarrow X\times \R$. Assume that $\Fix U \ne \emptyset$. Then $U$ is said to be
	\begin{enumerate}
	\item \emph{firmly nonexpansive} if
	\begin{equation*}
	\|U(\vec{x})-U(\vec{y})\|^2 + \|(\Id-U)(\vec{x})-(\Id-U)(\vec{y})\|^2 \leq \|\vec{x}-\vec{y}\|^2 \;\; \forall \vec{x} \in D, \;\; \forall \vec{y} \in D;
	\end{equation*}
	
	\item\emph{nonexpansive} if it is Lipschitz continuous with constant $1$,
	\begin{equation*}
	\|U(\vec{x})-U(\vec{y})\| \leq \|\vec{x}-\vec{y}\| \qquad \forall \vec{x} \in D, \quad \forall \vec{y} \in D;
	\end{equation*}
	
	\item\emph{quasinonexpansive (QNE)} if 
	$$
	\qquad \|U(\vec{x})-\vec{y}\| \leq \|\vec{x}-\vec{y}\| \qquad \forall \vec{x} \in D, \quad \forall \vec{y} \in \Fix U 
	$$ 
	(an operator that is both quasinonexpansive and continuous is called paracontracting);
	
	\item \emph{firmly quasinonexpansive (FQNE)} (or a \emph{cutter}) if 
	\begin{equation*}
	\|U\vec{x}-\vec{y}\|^2 + \|U\vec{x}-\vec{x}\|^2 \leq \|\vec{x}-\vec{y}\|^2\quad \forall \vec{x} \in D, \quad \forall \vec{y} \in \Fix U;
	\end{equation*}
	
	\item\emph{strictly quasinonexpansive (SQNE)} if
	\begin{equation*}
	\|U(\vec{x})-\vec{y}\| < \|\vec{x}-\vec{y}\| \qquad \forall \vec{x} \in D \setminus \Fix U,\quad \forall \vec{y} \in \Fix U;
	\end{equation*}
	\item\emph{$\rho$-strongly quasinonexpansive} for $\rho > 0$ if $$
	\|U\vec{x}-\vec{y}\|^2 \leq \|\vec{x}-\vec{y}\|^2 - \rho \|U\vec{x} -\vec{x}\|^2 \qquad \forall \vec{x} \in D \setminus \Fix U,\quad \forall \vec{y} \in \Fix U.
	$$
	\end{enumerate}
\end{definition}

\begin{lemma}{{\cite[Proposition~4.4]{BC}}}\label{lem:BCQNE}
	Let $D \subset X \times \R$ be nonempty. Let $U:D \rightarrow X\times \R$. The following are equivalent:
	\begin{enumerate}
		\item $U$ is firmly quasinonexpansive;
		\item $2U-\Id$ is quasinonexpansive;
		\item $(\forall \vec{x} \in D)(\forall \vec{y} \in \Fix U)\quad \|U\vec{x}-\vec{y}\|^2 \leq \langle \vec{x}-\vec{y}, U\vec{x}-\vec{y}\rangle$;
		\item $(\forall \vec{x} \in D)(\forall \vec{y} \in \Fix U)\quad \langle \vec{y}-U\vec{x}, \vec{x}-U\vec{x}\rangle \leq 0$;
		\item $(\forall \vec{x} \in D)(\forall \vec{y} \in \Fix U)\quad \|U\vec{x}-\vec{x}\|^2 \leq \langle \vec{y}-\vec{x},U\vec{x}-\vec{x}\rangle$.
	\end{enumerate}
\end{lemma}

\begin{definition}[Fej\'er monotonicity {\cite[5.1]{BC}}]A sequence $(\vec{x}_n)_{n\in \N}$ is Fej\'er monotone with respect to a closed convex set $C \subset X \times \R$ if
	$$
	\|\vec{x}_{n+1}-\vec{x}\|\leq \|\vec{x}_n - \vec{x}\|\qquad \forall \vec{x}\in C, \quad \forall n \in \N.
	$$
\end{definition}
A Fej\'er monotone sequence with respect to a closed convex set $C$ may be thought of as a sequence defined by $\vec{x}_n := U^n \vec{x}_0$ where $U$ is QNE with respect to $C=\Fix U$. We will make use of the fact that a Fej\'er monotone sequence with respect to a non-empty set is always bounded. We will also make use of the following convergence result.

\begin{theorem}{{\cite[Theorem 5.11]{BC}}}\label{thm:baucomstrongconvergence}
	Let $(\vec{x}_n)_{n\in\mathbb{N}}$ be a sequence in $X\times \R$ and let $C$ be a nonempty closed convex subset of $X\times \R$. Suppose that $(\vec{x}_n)_{n\in\mathbb{N}}$ is Fej\'er monotone with respect to $C$. Then the following are equivalent: 
	\begin{enumerate}
		\item the sequence $(\vec{x}_n)_{n\in\mathbb{N}}$ converges strongly (i.e. in norm) to a point in $C$;
		\item $(\vec{x}_n)_{n\in\mathbb{N}}$ possesses a strong sequential cluster point in $C$;
		\item $\underset{n\to \infty}{\rm \liminf}\; d(\vec{x}_n,C)=0.$
	\end{enumerate}
\end{theorem}

We will make use of the following result, which may be recognized as a simplified version of \cite[Theorem~4]{DLR18} and variants of which may be found in \cite{Cegielski}.

\begin{lemma}\label{lem:DLRcutterresult}
	Let $U: X\times R \rightarrow X \times \R$ be FQNE. Then the operator given by
	$$
	(2-\gamma)(U-\Id)+\Id
	$$
	is QNE for all $\gamma \in \left[0,2\right]$. Moreover, if $\gamma \in \left]0,2 \right[$ then $U$ is SQNE and the sequence $(\vec{x}_n)_{n \in \N}$ given by $\vec{x}_n := U^n \vec{x}_0$ satisfies
	\begin{equation*}
	\underset{n \rightarrow \infty}{\lim}\|\vec{x}_n - U(\vec{x}_n)\| \rightarrow 0.
	\end{equation*}
\end{lemma}

\section{The polar envelope and proximity operator}\label{sec:polarprox}

Friedlander, Mac\^{e}do, and Pong introduced the polar envelope and its associated polar proximity operator \cite{friedlander2019polar}, which we now recall.

\begin{definition}[Polar envelope and polar proximal map {\cite{friedlander2019polar}}]\label{def:pprox}
	
	For any closed gauge $\kappa: X\times \R \rightarrow \R$ and positive scalar $\alpha$, the function
	\begin{equation*}
	\kappa_\alpha: X \times \R \rightarrow \R: \vec{x} \mapsto \underset{\vec{u}}{\inf}\max \{\kappa(\vec{u}),(1/\alpha)\|\vec{x}-\vec{u}\|\}
	\end{equation*}
	is the \emph{polar envelope} of $\kappa$. The corresponding polar proximal map
	\begin{equation*}
	T_{\kappa, \alpha }: X\times \R \rightarrow X \times \R: \vec{x} \mapsto \underset{\vec{u}}{\argmin}\max\{\kappa(\vec{u}),(1/\alpha)\|\vec{x}-\vec{u}\| \}
	\end{equation*}
	sends a point $\vec{x}$ to the minimizing set that defines $\kappa_\alpha(\vec{x})$. Naturally, 
	\begin{equation*}
		f_\alpha^\pi: X \times \R \rightarrow \R: \vec{x} \mapsto \underset{\vec{u}}{\inf}\max \{\kappa(\vec{u}),(1/\alpha)\|\vec{x}-\vec{u}\|\}
	\end{equation*}
	denotes the polar envelope of the closed perspective transform $f^\pi$ for a proper convex function $f:X \rightarrow \R$.
\end{definition}

\begin{figure}[ht]
	\begin{center}
		\includegraphics[width=.8\textwidth]{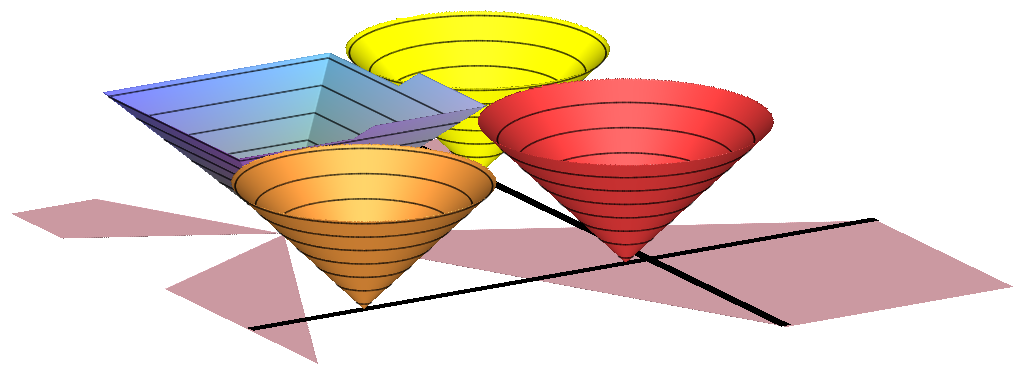}\\
		\includegraphics[width=.4\textwidth]{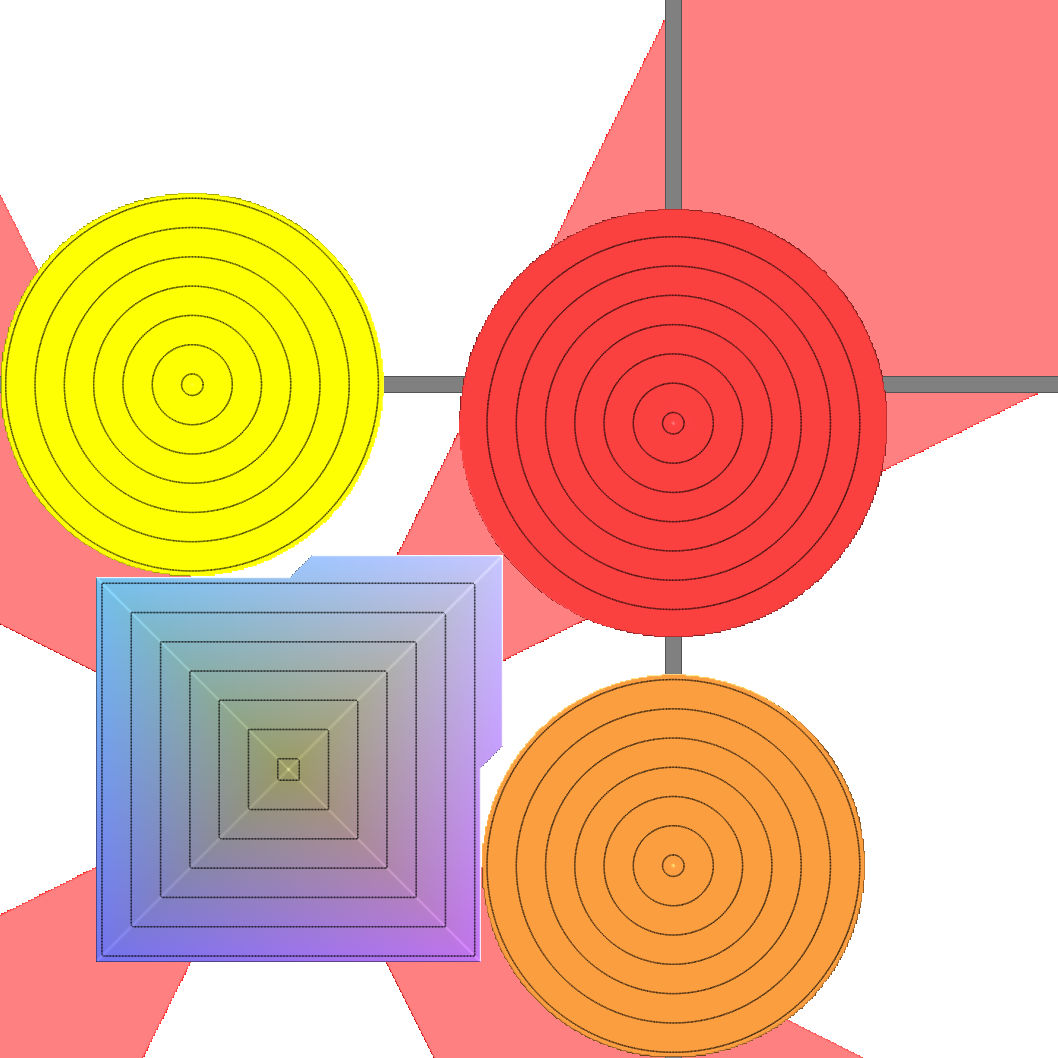}\includegraphics[width=.4\textwidth]{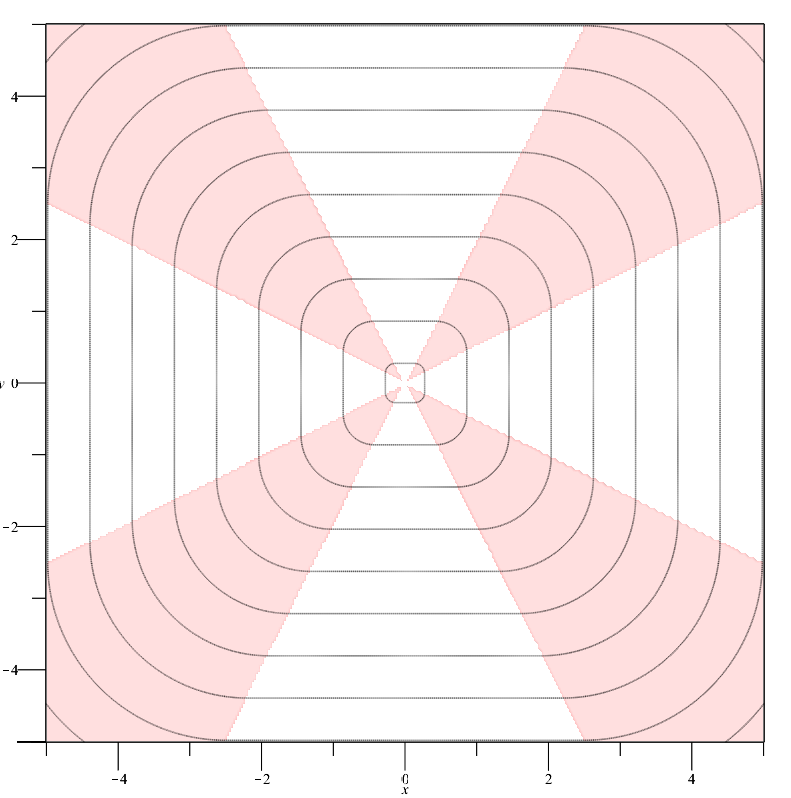}
	\end{center}
	\caption{Construction of the polar envelope and corresponding polar proximal map for $\kappa=\|\cdot\|_\infty$.}\label{fig:inftynorm}
\end{figure}

Figure~\ref{fig:inftynorm} shows the construction of the polar envelope and its proximity operator for $\kappa=\|\cdot \|_\infty$. At top and at bottom left, we take three choices of $\vec{x}$ and plot the functions $\|\cdot - \vec{x}\|$ in yellow, red, and orange respectively. The domain points for which each of these epigraphs intersects the epigraph of $\|\cdot \|_\infty$ at lowest height are the respective proximal points. The height at the point of intersection determines the envelope value. For points $\vec{x}$ in the white regions, such as the points for which the functions $\|\cdot-\vec{x}\|$ are orange and yellow respectively, the envelope value is simply $\|\vec{x}\|_\infty/2$. For points lying in the red regions, the proximal point lies on the diagonals; for points in the interiors of the red regions, the envelope values are strictly greater than $\|\vec{x}\|_\infty/2$. This results in the smoothing apparent in the red regions for the envelope shown at bottom right.

We devote the remainder of this section to showing that the polar proximity operator $T_{\kappa,\alpha}$ is firmly quasinonexpansive. We have the symmetry of vertical rescaling,
\begin{align*}
T_{\kappa,\alpha}(\vec{u})&=\argminn{\vec{u}} \max \left\{\kappa(\vec{u}),(1/\alpha)\|\vec{x}-\vec{u}\| \right\}\\
&=\argminn{\vec{u}} \max \left\{\alpha \kappa(\vec{u}),\|\vec{x}-\vec{u}\| \right\}=T_{\alpha\kappa,1}(\vec{u}),
\end{align*}
by which the proximity operators satisfy $T_{\kappa,\alpha} = T_{\alpha \kappa,1}$, while the envelopes satisfy $(\alpha \kappa)_1 = \alpha(\kappa_\alpha)$. Thus, by working with a general $\kappa$, we can, and do, let $\alpha=1$ without loss of generality. For simplicity, we also write $T$ instead of $T_{\kappa,\alpha}$. Note that we still need the notation $\kappa_\alpha$ to distinguish the polar envelope from the gauge $\kappa$ itself.

\begin{lemma}\label{lem:FixT}
	For a closed gauge $\kappa$, it holds that
	\begin{equation*}%\label{FixT}
	\Fix T = \zer \kappa.
	\end{equation*}
\end{lemma}
\begin{proof}
	The fact that $\zer \kappa \subset \Fix T$ is obvious. We will show the reverse inclusion. Let $\vec{x} \in \Fix T$. Then 
	\begin{equation*}
	\vec{x} = T\vec{x} = \underset{\vec{u}\in X\times \R}{\argmin}\max \{\kappa(\vec{u}) , \|\vec{x}-\vec{u}\|\},
	\end{equation*}
	and so 
	\begin{equation*}
	\underset{\vec{u} \in X \times \R}{\inf}\max \{\kappa(\vec{u}) , \|\vec{x}-\vec{u}\|\} = \max\{\kappa(\vec{x}) , \|\vec{x}-\vec{x}\|\} = 0.
	\end{equation*}
	Thus $\kappa(\vec{x}) \leq 0$. Combining with the fact $\kappa(\vec{x})\geq 0$, we have $\kappa(\vec{x}) = 0$.
\end{proof}

\begin{lemma}\label{lem:kappazero}
	Let $\kappa$ be a closed gauge and $\kappa_\alpha$ its polar envelope. Then
	$$(\kappa(\vec{x})=0) \iff (\kappa_\alpha(\vec{x}) = 0).$$
\end{lemma}
\begin{proof}
	Let $\kappa(\vec{x}) = 0$. Then 
	$$
	\kappa_\alpha(\vec{x})= \underset{\vec{u} \in X \times \R}{\inf}\max \{\kappa(\vec{u}),\|\vec{x}-\vec{u}\| \} \leq \max \{ \kappa(\vec{x}),\|\vec{x}-\vec{x}\|\} = 0.
	$$
	Thus $\kappa_\alpha(\vec{x}) = 0$. 
	
	Now let $\kappa_\alpha(\vec{x}) = 0$. Since $\kappa_\alpha(\vec{x}) = 0$, there exists a sequence $(\vec{u}_n)_{n \in N}$ such that 
	$$
	\max \{ \kappa(\vec{u}_n), (1/\alpha)\|\vec{x}-\vec{u}_n\| \} \rightarrow 0.
	$$
	Then we have that $\kappa(\vec{u}_n) \rightarrow 0$ and $\|\vec{x}-\vec{u}_n\| \rightarrow 0$. Since $\|\vec{x}-\vec{u}_n\| \rightarrow 0$, we have that $\vec{u}_n \rightarrow x$. Combining with the fact that $\kappa$ is lower semicontinuous, we have that $\kappa(\vec{x}) \leq \lim_{\vec{u}_n \rightarrow \vec{x}} \kappa(\vec{u}_n) = 0$. This concludes the result.
\end{proof}

%\begin{remark}[{\cite[Section~5]{friedlander2019polar}}]\label{rem:subdifferentialkappa}
%	Let $\kappa$ be a closed gauge and $\kappa_\alpha(x) > 0$. Then one of the following holds:
%	\begin{enumerate}[label=(\roman*)]
%		\item $\kappa(Tx) <\frac{1}{\alpha}\|Tx-x\|$, in which case $0 \in \frac{1}{\alpha} \frac{Tx-x}{\|Tx-x\|} + N_{\dom \kappa}(Tx)$ and $Tx = P_{\dom \kappa}x$.
%		\item We have that
%		\begin{align*}
%		r &:= \kappa(Tx) = \frac{1}{\alpha}\|Tx-x\|\\
%		\text{where}\quad Tx &= P_{\lev_{\leq r}(\kappa)}x, 
%		\end{align*} 
%		and there exists $\lambda \in \left [0,1 \right[$ such that 
%		\begin{align*}
%		0 &\in \frac{1-\lambda}{\alpha} \frac{Tx-X}{\|Tx-x\|} + \partial (\lambda^+ \kappa)(Tx) \nonumber \\
%		& = \frac{1-\lambda}{\alpha} \frac{Tx-x}{\|Tx-x\|} + \begin{cases}
%		\left \{ \lambda z \; | \; z \in \partial \kappa (Tx) \right  \} & \text{if}\; \lambda > 0 \\
%		N_{\dom \kappa}(Tx) & \text{if}\; \lambda = 0
%		\end{cases}
%	\end{align*}
%	\end{enumerate}
%\end{remark}

\begin{lemma}\label{lem:subdifferentialkappa}
	Let $\kappa$ be a closed gauge and $\kappa(\vec{x}) > 0$. Then one of the following holds:
	\begin{enumerate}[label=(\roman*)]
		\item\label{case1} $\kappa(T\vec{x}) <\frac{1}{\alpha}\|T\vec{x}-\vec{x}\|$, in which case $0 \in \frac{1}{\alpha} \frac{T\vec{x}-\vec{x}}{\|T\vec{x}-\vec{x}\|} + N_{\dom \kappa}(T\vec{x})$ and $T\vec{x} = P_{\dom \kappa}\vec{x}$.
		\item\label{case2} We have that
		\begin{align*}
		r :=& \kappa(T\vec{x}) = \frac{1}{\alpha}\|T\vec{x}-\vec{x}\|\\
		\text{where}\quad T\vec{x} =& P_{\lev_{\leq r}(\kappa)}\vec{x}, 
		\end{align*}
		and there exists $\lambda \in \left [0,1 \right[$ such that 
		\begin{align*}
		0 &\in \frac{1-\lambda}{\alpha} \frac{T\vec{x}-\vec{x}}{\|T\vec{x}-\vec{x}\|} + \partial (\lambda^+ \kappa)(T\vec{x}) \nonumber \\
		& = \frac{1-\lambda}{\alpha} \frac{T\vec{x}-\vec{x}}{\|T\vec{x}-\vec{x}\|} + \begin{cases}
		\left \{ \lambda z \; | \; z \in \partial \kappa (T\vec{x}) \right  \} & \text{if}\; \lambda > 0 \\
		N_{\dom \kappa}(T\vec{x}) & \text{if}\; \lambda = 0
		\end{cases}.
		\end{align*}
	\end{enumerate}
\end{lemma}
\begin{proof}
	Let $\kappa(\vec{x})>0$. From Lemma~\ref{lem:kappazero}, $\kappa_\alpha(\vec{x})>0$. We have from \cite[Section~5]{friedlander2019polar} that the condition $\kappa_\alpha(\vec{x})>0$ guarantees that one of \ref{case1} or \ref{case2} must hold.
\end{proof}

We will use the characterization of $T$ from Lemma~\ref{lem:subdifferentialkappa} often; hence the shorthand $P_{\kappa \leq r}:=P_{\lev_{\leq r}(\kappa)}$. Now we have the principle result of this section, which establishes that the polar proximity operator $T$ is firmly quasinonexpansive.

\begin{theorem}\label{thm:TQNE}$T$ is firmly quasinonexpansive.
\end{theorem}
\begin{proof}
	Let $\vec{y} \in X \times \R$ and $\vec{x} \in \Fix T$. If $\kappa(\vec{y}) = 0$ then $\vec{y} \in \Fix T$, so let $\kappa(\vec{y}) > 0$. By Lemma~\ref{lem:subdifferentialkappa}, we need only consider two cases. 
	
	\textbf{Case 1}: If Lemma~\ref{lem:subdifferentialkappa}\ref{case1} holds, then we have that
	\begin{equation}\label{case:normalcone}
	-\frac{1}{\alpha}\frac{T\vec{y}-\vec{y}}{\|T\vec{y}-\vec{y}\|} \in N_{\dom \kappa} (T\vec{y}),
	\end{equation}
	and so $\vec{y} - T\vec{y} \in N_{\dom \kappa}(T\vec{y})$. Thus $T\vec{y} = P_{\cdk}(\vec{y})$. The operator $P_{\cdk}$ is FQNE with $\Fix P_{\cdk} = \cdk$. Thus by Lemma~\ref{lem:BCQNE},
	\begin{equation}\label{projFQNE}
	(\forall \vec{v} \in X \times \R)\left(\forall \vec{u} \in \Fix P_{\cdk}\right)\quad \left \langle \vec{u}-P_{\cdk}\vec{v},\vec{v}-P_{\cdk}\vec{v} \right\rangle \leq 0.
	\end{equation}
	Since $\vec{x} \in \Fix T \subset \cdk = \Fix P_{\cdk}$, we have that \eqref{projFQNE} is true, in particular, for $\vec{u}=\vec{x}$ and $\vec{v}=\vec{y}$. Thus we obtain
	$$
	\left \langle \vec{x}-P_{\cdk}\vec{y},\vec{y}-P_{\cdk}\vec{y} \right \rangle \leq 0.
	$$
	This is just 
	$$
	\langle \vec{x}-T\vec{y},\vec{y}-T\vec{y} \rangle \leq 0.
	$$
	By Lemma~\ref{lem:BCQNE}, this is what we needed to show.
	
	\textbf{Case 2}: If Lemma~\ref{lem:subdifferentialkappa}\ref{case2} holds with $\lambda = 0$, then we again obtain \eqref{case:normalcone} and proceed as in Case 1, obtaining what we needed to show. If Lemma~\ref{lem:subdifferentialkappa}\ref{case2} holds with $\lambda > 0$ then there exists $r >0$ such that
	\begin{equation*}%\label{case:projlev}
	T\vec{y} = P_{\kappa \leq r}\vec{y}.
	\end{equation*}
	Now since $P_{\kappa \leq r}$ is FQNE, we have that 
	\begin{equation}\label{projlevFQNE}
	(\forall \vec{v} \in X \times \R)\left(\forall \vec{u} \in \Fix P_{\kappa \leq r}\right)\quad  \left\langle \vec{u}-P_{\kappa \leq r}\vec{v} , \vec{v}-P_{\kappa \leq r}\vec{v} \right\rangle \leq 0.
	\end{equation}
	Since $\vec{x} \in \Fix T =\lev_{\leq0}\kappa \subset \lev_{\leq r_2}\kappa = \Fix P_{\kappa \leq r_2}$, we have that \eqref{projlevFQNE} is true, in particular, for $\vec{u}=\vec{x}$ and $\vec{v}=\vec{y}$. Thus we obtain
	$$
	\langle \vec{x}-P_{\kappa \leq r_2}\vec{y},\vec{y}-P_{\kappa \leq r_2}\vec{y} \rangle \leq 0.
	$$
	This is just 
	$$
	\langle \vec{x}-T\vec{y},\vec{y}-T\vec{y} \rangle \leq 0.
	$$
	By Lemma~\ref{lem:BCQNE}, this is what we needed to show.
\end{proof}

\section{The projected polar proximal point algorithm}\label{sec:PPPPA}

We now recall the projected polar proximal point algorithm.

\begin{definition}[Projected polar proximal point algorithm \PPPPA \; {\cite[7.2]{friedlander2019polar}}]\label{def:p4}
	Fix $\alpha>0$ and set
	\begin{align*}
	%\rho_{\alpha,f}: X \rightarrow X : v \mapsto &f_\alpha^\pi(v,1),\\
	\mathfrak{P}_{\alpha,f}: X \rightarrow X: v \mapsto&T_{f^\pi,\alpha}(v,1).
	\end{align*}
	The projected polar proximal point algorithm is to begin with any $v_0$ and update by
	\begin{equation}\label{seq:PPPPA}
	(v_{k+1},\lambda_{k+1})=\mathfrak{P}_{\alpha,f}(v_k).
	\end{equation}
\end{definition}
Intuitively, the motivation of \PPPPA \;is to minimize a function $f$ by attacking the gauge given by its closed perspective transform $f^\pi$. Its polar envelope $f_\alpha^\pi$ serves a role analogous to the role played by the Fenchel--Moreau envelope in the construction of the traditional proximal point algorithm; see the remarks in \cite[7.2]{friedlander2019polar}. The use of the gauge allows problems to be reformulated using gauge duality \cite{friedlander2014gauge}.

In addition to explaining these connections, Friedlander, Mac\^{e}do, and Pong also showed that \PPPPA \;has a useful fixed point property, which we now recall.

\begin{theorem}[Fixed points of \PPPPA \; {\cite[Theorem~7.4]{friedlander2019polar}}]\label{thm:fixedpoints_argmin}
	Let $f:X \rightarrow \R_+ \cup \{+\infty \}$ be a proper closed nonnegative convex function with $\inf f>0$ and $\argmin f \neq \emptyset$. The following hold
	\begin{enumerate}[label=(\roman*)]
		\item\label{fixedpoints1} If $(v,\lambda_*)=\mathfrak{P}_{1,f}$, then $\lambda_* > 0$ and $\lambda_*^{-1}v \in \argmin f$.
		\item\label{fixedpoints2} If $v \in \argmin f$, then there exists $\lambda_* > 0$ so that $(\tau v, \lambda_*) = \mathfrak{P}_{1,f}(\tau v)$ where $\tau := \left[1+f(v) \right]^{-1}$.
	\end{enumerate}
\end{theorem}

To show convergence of \PPPPA, we will analyse a generalization of it, which we now introduce.

\begin{definition}[Generalized projected polar proximal point algorithm (\GPPPPA)]\label{def:GPPPPA}
	For a gauge $\kappa: X\times \R \rightarrow \R$ and fixed $\alpha>0$, choose a starting point $\vec{x}_0 \in X\times \R$ and iterate by
	\begin{align}
	\vec{x}_{k+1} &= P_S \circ T \vec{x}_k,\label{seq:GPPPPA}\\
	\text{where}\quad S &= X \times \{1\}, \nonumber\\
	\text{and}\quad P_S&: X\times \R \rightarrow X \times \R: (y,\lambda) \mapsto (y,1) \nonumber
	\end{align}
	is simply the projection operator for $S$.
\end{definition}
Note that \PPPPA \;and \GPPPPA \;are defined on $X$ and $X \times \R$ respectively. However, when $\kappa = f^\pi$ the sequences $(v_k)$ from \eqref{seq:PPPPA} and $(\vec{x}_k)_k$ from \eqref{seq:GPPPPA} clearly satisfy $(v_k,1) = \vec{x}_k$, and so the algorithms \PPPPA \;and \GPPPPA \;generate the same sequence on the non-lifted space. This means that we can study \PPPPA \;by studying \GPPPPA, because their performance for $\kappa = f^\pi$ is the same.

For our purposes, the advantage of \GPPPPA \;is that it is defined on the lifted space $X \times \R$. This allows us to decompose the method into iterative application of the two separate operators: $T$ and $P_S$. This allows for greater flexibility. Specifically, it allows us to study \GPPPPA \;as a splitting method, where two different operators are applied in succession: first the one and then the other. This allows us to build new characterizations of fixed points, and these characterizations, in turn, allow us to show global convergence. Thus, by analysing \GPPPPA, we are able to prove the desired convergence of \PPPPA\; without assuming strong convexity of $(f^\pi)^2$. Whether the added flexibility of \GPPPPA \;has benefits beyond its utility for learning about \PPPPA \;is a natural question for future research. For the present, our main motivation for introducing \GPPPPA \;is the aforementioned advantage.

\subsection{Alternative Fixed Point Characterization}\label{sec:fixedpoints}

We next establish a useful characterization of the fixed points (Proposition~\ref{prop:fixedpoints}). For the purpose, we need the following two lemmas.

\begin{lemma}\label{lem:abslambdaleq1}
	Let $(x,1) \in \Fix P_S \circ T$. Then the following hold.
	\begin{enumerate}[label=(\roman*)]
		\item\label{eqn:abslambdaleq1} $T(x,1) = (x,\lambda)$ for some $\lambda \in \left[0,1 \right]$.
		\item\label{eqn:abslambdaless1} Moreover, if $(x,1) \notin \Fix T$, then $\lambda \in \left[0,1\right[$.
	\end{enumerate}
\end{lemma}
\begin{proof}
	\ref{eqn:abslambdaleq1}: Since $(x,1) \in \Fix P_S \circ T$, we have that 
	$$
	P_S \circ T(x,1) = (x,1).
	$$
	Since $P_S^{-1}(w,1) = \{(w,\lambda)\; | \lambda \in \R \}$ for all $w \in X$, we have that $T(x,1) = (x,\lambda)$ for some $\lambda \in \R$. We now show that $\lambda \in \left[0,1\right]$. Since $T$ is $FQNE$,
	\begin{equation}\label{eqn:TQNE}
	(\forall \vec{u} \in X\times \R)(\forall \vec{v}\in \Fix T)\quad \langle \vec{v}-T\vec{u},\vec{u}-T\vec{u} \rangle \leq 0. 
	\end{equation}
	In particular, $\vec{0} \in \Fix T$ and so \eqref{eqn:TQNE} yields
	\begin{equation}\label{eqn:normcontracting}
	\langle -T\vec{u},\vec{u}-T\vec{u} \rangle \leq 0. 
	\end{equation}
	In particular \eqref{eqn:normcontracting} holds for $\vec{u}=(x,1)$, and so we have
	\begin{align*}
	\langle -T(x,1) , (x,1) - T(x,1) \rangle \leq 0
	\end{align*}
	This is just
	\begin{align*}
	\langle -(x,\lambda),(x,1)-(x,\lambda) \rangle &= \langle -(x,\lambda),(0,1-\lambda)\rangle \\
	&= -\lambda(1-\lambda) \leq 0
	\end{align*}
	Thus $\lambda \in \left[0,1 \right]$.
	
	\ref{eqn:abslambdaless1}: Since $T$ is FQNE, it is SQNE \cite{BC}, and so we have that
	\begin{equation}\label{eqn:normcontractinghard}
	(\forall \vec{u} \in X\times \R \setminus \Fix T)(\forall \vec{v} \in \Fix T)\quad \|T\vec{u}-\vec{v}\|^2 < \|\vec{u}-\vec{v}\|^2.
	\end{equation}
	Specifically, \eqref{eqn:normcontractinghard} holds for $\vec{u} = (x,1)$ and $\vec{v}=\vec{0}$, and so we obtain
	$$
	\|(x,\lambda)\|^2 < \|(x,1)\|^2.
	$$ 
	This is just
	$$
	\|x\|^2 + \lambda^2 < \|x\|^2 + 1^2,
	$$
	and so we conclude that $\lambda < 1$. This concludes the result.
\end{proof}

\begin{lemma}\label{lem:fgreaterthanb}Let $(x,1) \in \Fix P_S \circ T$ and $(x,\lambda) := T(x,1)$ (a representation that always holds by Lemma~\ref{lem:abslambdaleq1}). Then for any $(y,1) \in S$ the following hold:
	\begin{enumerate}[label=(\roman*)]
		\item\label{lem:fgreaterthanbitem1} $\|T(y,1)-(y,1)\| \geq \|T(x,1)-(x,1) \|$;
		\item\label{lem:fgreaterthanbitem2} Additionally, if $T(x,1) = P_{\cdk}(x,1)$ or $\lambda=1$ then 
		\begin{enumerate}
			\item\label{lem:fgreaterthanbitem2a} $\|T(y,1)-P_S\circ T(y,1)\| \geq \|T(x,1)-(x,1)\|$;
			\item\label{lem:fgreaterthanbitem2b} $T(y,1) = (w,\mu)$ for some $\mu \in \R$ satisfying $|1-\lambda| \leq |1-\mu|$;
			\item\label{lem:fgreaterthanbitem2c}If $\lambda < 1$, then $T(y,1) = (w,\mu)$ for some $\mu \leq \lambda$. 
		\end{enumerate}	
	\end{enumerate}
\end{lemma}
\begin{proof}
	Let $(y,1) \in S$ and $(w,\mu) = T(y,1)$. Set $r':=\|T(x,1)-(x,1) \|$ and $r_2:=\|T(y,1)-(y,1)\|$. We have by Lemma~\ref{lem:abslambdaleq1} that 
	\begin{equation*}%\label{Pr1cutter3}
	\lambda \in \left[0,1 \right].
	\end{equation*}
	We will consider two cases: when $\lambda = 1$ and when $\lambda \in \left[0,1\right[$.
	
	\textbf{Case $\lambda =1$}. Then $|1-\lambda| = 0$, and so \ref{lem:fgreaterthanbitem2b} clearly holds. Additionally, by Lemma~\ref{lem:abslambdaleq1}, we have that $(x,1) \in \Fix T$, and so $\|T(x,1)-(x,1)\| = 0$, and so  \ref{lem:fgreaterthanbitem1} and \ref{lem:fgreaterthanbitem2a} both clearly hold, while \ref{lem:fgreaterthanbitem2c} does not apply. This concludes what we needed to show in the case $\lambda=1$.
	
	\textbf{Case $\lambda <1$}. Let
	\begin{equation}\label{lambdalessthanone}
	\lambda \in \left[0,1 \right[
	\end{equation}
		
	By Lemma~\ref{lem:subdifferentialkappa}, we may further consider two subcases, namely when  $T(x,1) = P_{\cdk}(x,1)$ and when $T(x,1) = P_{\kappa \leq r'}(x,1)$.
	
	\textbf{Subcase 1}: Let $T(x,1) = P_{\cdk}(x,1)$. Since $P_{\cdk}$ is FQNE with $\Fix P_{\cdk} = \cdk$, we have that
	\begin{equation}\label{projdomcutter}
	(\forall \vec{u} \in X\times \R)(\forall \vec{v} \in \cdk) \quad \langle P_{\cdk}\vec{u}-\vec{u},P_{\cdk}\vec{u}-\vec{v} \rangle \leq 0.
	\end{equation}
	In particular, \eqref{projdomcutter} holds for $\vec{u} = (x,1)$ and $\vec{v} = T(y,1)$. Thus
	\begin{equation}\label{projdomcutter2}
	\langle P_{\cdk}(x,1)-(x,1) , P_{\cdk}(x,1) - T(y,1)\rangle \leq 0.
	\end{equation}
	
	Since $P_{\cdk}(x,1) = T(x,1)$, we have $P_{\cdk}(x,1) = (x,\lambda)$. Combining with \eqref{projdomcutter2}, we obtain
	\begin{align}
	\langle (x,\lambda)-(x,1) , (x,\lambda)-T(y,1) \rangle &\leq 0 \label{quadraticinlambdatop}\\
	\langle (0,\lambda-1) , (x,\lambda)-T(y,1) \rangle &\leq 0\nonumber \\
	(\lambda-1)(\lambda-\mu) &\leq 0.\label{quadraticinlambda}
	\end{align}
	Recall that by \eqref{lambdalessthanone}, we have $\lambda \in \left[0 ,1\right[$. Thus we have $(\lambda -1) < 0$. Combining this fact with \eqref{quadraticinlambda}, we have $\lambda - \mu \geq 0$, and so $\lambda \geq \mu$. This shows \ref{lem:fgreaterthanbitem2b} and \ref{lem:fgreaterthanbitem2c}. Thus $(1-\mu) \geq (1-\lambda) > 0$, where the final inequality is because $\lambda < 1$. Thus we have 
	\begin{equation}\label{biginequality1}
	(\mu-1)^2 \geq (\lambda - 1)^2. 
	\end{equation}
	We also have that 
	\begin{equation}\label{biginequality2}
	\|T(y,1)-(y,1)\|^2 = \|(w,\mu)-(y,1) \|^2 = \|w-y\|^2+(\mu-1)^2 \geq (\mu-1)^2, 
	\end{equation}
	and that
	\begin{equation}\label{biginequality3}
	(\lambda - 1)^2 = \|(x,1)-(x,\lambda)\|^2 = \|T(x,1)-(x,1)\|^2 = r'^2.
	\end{equation}
	Combining \eqref{biginequality1}, \eqref{biginequality2}, and \eqref{biginequality3}, we obtain $\|T(y,1) - (y,1)\| \geq r'$. Thus $r_2 \geq r'$, which shows \ref{lem:fgreaterthanbitem1}.	
	
	We also have that
	\begin{align*}
	\|T(x,1)-P\circ T(x,1)\| &= \|(x,\lambda)-(x,1)\| = 1-\lambda \leq 1-\mu\\
	\text{and}\quad 1-\mu &= \|(w,\mu)-(w,1)\|= \|T(y,1)-P\circ T(y,1)\|,
	\end{align*}
	which shows \ref{lem:fgreaterthanbitem2a}. Thus we have shown everything we needed to show in the case when $T(x,1) = P_{\cdk}(x,1)$

	\textbf{Subcase 2}: Let $T(x,1) = P_{\kappa \leq r'}(x,1)$. Suppose for a contradiction that $r_2 < r'$. For the sake of simplicity, define $P_{r'}:=P_{\kappa \leq r'}$. We have that
	$$
	\kappa(T(y,1)) \leq r_2 < r',
	$$
	where the first inequality is because $\kappa(T(y,1)) \leq \|T(y,1)-(y,1)\|$ \cite[Theorem~4.4]{friedlander2019polar} and the second is our contradiction assumption. Thus 
	$$
	T(y,1) \in \lev_{\kappa \leq r_2} \kappa \subset \lev_{\kappa < r'} \subset \lev_{\kappa \leq r'}.$$ 
	
	As $\lev_{\leq r'}\kappa$ is closed and convex, the operator $P_{r'}$ is FQNE with $\Fix P_{r'} = \lev_{\leq r'}\kappa$. Using Lemma~\ref{lem:BCQNE}, we have that
	\begin{equation}\label{Pr1cutter}
	(\forall \vec{u} \in X \times \R)(\forall \vec{v} \in \lev_{\leq r'}\kappa )\quad\langle P_{r'}\vec{u}-\vec{u} , P_{r'}\vec{u}-\vec{v} \rangle \leq 0.
	\end{equation}
	In particular, $T(y,1) \in \Fix P_{\kappa \leq r'}$, and so we can apply \eqref{Pr1cutter} with $\vec{u} = (x,1)$ and $\vec{v} = T(y,1)$, obtaining
	\begin{equation}\label{Pr1cutter2}
	\langle P_{r'}(x,1)-(x,1) , P_{r'}(x,1) - T(y,1)\rangle \leq 0.
	\end{equation}
	Since $P_{r'}(x,1) = T(x,1)$, we have $P_{r'}(x,1) = (x,\lambda)$. Let $(w,\mu):= T(y,1)$. Combining with \eqref{Pr1cutter2}, we again obtain \eqref{quadraticinlambdatop}, and we proceed as we did in Case 1 to obtain $r_2 \geq r'$, a contradiction. Thus $r_2 \geq r'$, which shows \ref{lem:fgreaterthanbitem1}.
\end{proof}

We now establish the useful alternative characterization of the set $\Fix P_S \circ T$. 

\begin{proposition}[Fixed points of $P_S\circ T$]\label{prop:fixedpoints}
	Let $(x,1) \in \Fix P_S \circ T$ and set $r' := \|(x,1)-T(x,1)\|$. It holds that 
	\begin{equation}\label{characterizationofr1}
	\Fix P_S \circ T = \{(u,1) \;|\; \|T(u,1)-(u,1)\| = r' \}.
	\end{equation}
\end{proposition}
\begin{proof}
	The first inclusion $\{(u,1) \; | \; \|T(u,1)-(u,1)\| = r' \} \supset \Fix P_S \circ T$ is a consequence of Lemma~\ref{lem:fgreaterthanb}. Simply let $(x,1),(y,1) \in \Fix P_S \circ T$, and we have from Lemma~\ref{lem:fgreaterthanb}\ref{lem:fgreaterthanbitem1} that 
	\begin{align*}\|T(y,1)-(y,1)\| &\geq \|T(x,1)-(x,1) \|,\\
	\text{and}\quad \|T(y,1)-(y,1)\| &\leq \|T(x,1)-(x,1)\|,
	\end{align*}
	and so $\|T(y,1)-(y,1)\| = \|T(x,1)-(x,1)\| = r'$. 
	
	Now we will show the reverse inclusion. Let $(y,1) \in S$ and let
	\begin{equation*}%\label{def:r1}
	\|(y,1)-T(y,1)\| = r' = \|(x,1)-T(x,1)\|.
	\end{equation*}
	By Lemma~\ref{lem:abslambdaleq1} we have that 
	\begin{equation*}
	\exists \lambda \in \left[0,1\right] \quad \text{such that}\quad (x,\lambda)=T(x,1).
	\end{equation*}
	
	\textbf{Case 1: $\lambda=1$}. Suppose $\lambda =1$. Then we have that $(x,1) \in \Fix T$ and so $r' = 0$. Since $\|T(y,1)-(y,1)\|=r'$, we have $\|T(y,1)-(y,1)\| = 0$. Thus $(y,1) \in \Fix T$, and so $(y,1) \in \Fix P_S \circ T$. This concludes the case when $\lambda=1$.
	
	\textbf{Case 2: $\lambda \in \left[0,1\right[$}. Let $\lambda \in \left[0,1\right[$. Then $(x,1) \notin \Fix T$ and so $\kappa(x,1) > 0$. Since $\kappa(x,1)>0$, we have by Lemma~\ref{lem:subdifferentialkappa} that either $T(x,1) = P_{\cdk}(x,1)$ or $T(x,1) = P_{\kappa \leq r'}(x,1)$.
	
	\textbf{Case 2(a): $T(x,1) = P_{\cdk}(x)$}. Since $\cdk$ is closed and convex, we have that $P_{\cdk}$ is FQNE. Since $P_{\cdk}$ is FQNE with $\Fix P_{\cdk} = \cdk$, we have that
	\begin{equation}\label{fixprojdomcutter}
	(\forall \vec{v} \in X\times \R)(\forall \vec{u} \in \cdk) \quad \langle \vec{u}-P_{\cdk}\vec{v},\vec{v}-P_{\cdk}\vec{v} \rangle \leq 0.
	\end{equation}
	In particular, $T(y,1) \in \cdk$, and so we may apply \eqref{fixprojdomcutter} holds with $\vec{u} = T(y,1)$ and $\vec{v} = (x,1)$, obtaining
	\begin{equation}\label{fixprojdomcutter2}
	\langle T(y,1)-P_{\cdk}(x,1), (x,1)-P_{\cdk}(x,1) \rangle \leq 0.
	\end{equation}
	Using the fact that $T(y,1) = (w,\mu)$ and $P_{\cdk}(x,1) = T(x,1) = (x,\lambda)$, \eqref{fixprojdomcutter2} becomes
	\begin{equation}\label{fixprojdomcutter3}
	\langle (w,\mu)-(x,\lambda),(x,1)-(x,\lambda) \rangle \leq 0.
	\end{equation}
	From \eqref{fixprojdomcutter3} we have that
	\begin{equation}\label{fixprojdomcutter4}
	(\mu - \lambda)(1-\lambda) \leq 0.
	\end{equation}
	Using the fact that $\lambda < 1$, \eqref{fixprojdomcutter4} implies that $\mu \leq \lambda$. Thus we have that $1-\mu \geq 1-\lambda \geq 0$. Thus we have that
	\begin{equation}\label{fixprojdomcutter5}
	(1-\mu)^2 \geq (1-\lambda)^2.
	\end{equation}
	Now since 
	\begin{align*}
	(r')^2 &= \|(x,1)-T(x,1)\|^2 = \|(x,1)-(x,\lambda)\|^2 = (1-\lambda)^2\quad \text{and}\\
	(r')^2 &= \|(y,1)-T(y,1)\|^2 = \|(y,1)-(w,\mu)\|^2 = \|y-w\|^2 + (1-\mu)^2,
	\end{align*}
	we have that
	\begin{equation}\label{fixprojdomcutter6}
	(1-\lambda)^2 = \|y-w\|^2 + (1-\mu)^2.
	\end{equation}
	Combining \eqref{fixprojdomcutter5} and \eqref{fixprojdomcutter6}, we obtain
	\begin{equation*}
	0\geq (1-\lambda)^2 - (1-\mu)^2 = \|w-y\|^2.
	\end{equation*}
	Thus we have that $w=y$, and so $T(y,1) = (y,\mu)$. Thus
	$$P_S\circ T(y,1) = P_S (y,\mu) = (y,1),$$
	and so $(y,1) \in \Fix P_S \circ T$.
	
	\textbf{Case 2(b): $T(x,1) = P_{\kappa\leq r'}(x,1)$}. Let $T(x,1) = P_{\kappa\leq r'}(x,1)$. 
	
	Since $P_{\kappa \leq r'}$ is FQNE with $\Fix P_{\kappa \leq r'} = \lev_{\leq r'}\kappa$, we have that
	\begin{equation}\label{fixprojlevcutter}
	(\forall \vec{v} \in X\times \R)(\forall \vec{u} \in \lev_{\leq r'} \kappa) \quad \langle \vec{u}-P_{\kappa \leq r'}\vec{v},\vec{v}-P_{\kappa \leq r'}\vec{v} \rangle \leq 0.
	\end{equation}
	In particular, since $\kappa(T(y,1))\leq \|T(y,1)-(y,1)\|=r'$, we have that $T(y,1) \in \Fix P_{\kappa \leq r'}$, and so we may apply \eqref{fixprojlevcutter} with $\vec{u} = T(y,1)$ and $\vec{v} = (x,1)$, obtaining
	\begin{equation}\label{fixprojlevcutter2}
	\langle T(y,1)-P_{\kappa \leq r'}(x,1), (x,1)-P_{\kappa \leq r'}(x,1) \rangle \leq 0.
	\end{equation}
	Using the fact that $T(y,1) = (w,\mu)$ and $P_{\kappa \leq r'}(x,1) = T(x,1) = (x,\lambda)$, \eqref{fixprojlevcutter2} again yields \eqref{fixprojdomcutter3}, and we proceed as in \textbf{Case 2(a)}. 
	
	This shows the desired result.
\end{proof}	
The following shorthand will simplify notation in the results that follow. Whenever $\Fix P_S \circ T \neq \emptyset$, we define
	\begin{equation*}\label{def:E}
	E: (\Fix P_S \circ T) \times S \rightarrow \R: \quad (\vec{x},\vec{y}) \mapsto \|T\vec{y}-\vec{y}\|^2 - \|T\vec{x}-\vec{x}\|^2.
	\end{equation*}

Our previous results admit the following important property.

\begin{lemma}\label{lem:Egeq0}Whenever $\Fix P_S \circ T \neq \emptyset$, we have that
	\begin{equation*}
	(\forall \vec{x} \in \Fix P_S \circ T)(\forall \vec{y} \in S)\quad E(\vec{x},\vec{y}) \geq 0,
	\end{equation*}
	and equality holds if and only if $\vec{y} \in \Fix P_S \circ T$.
\end{lemma}
\begin{proof}
	The inequality is an immediate consequence of Lemma~\ref{lem:fgreaterthanb}\ref{lem:fgreaterthanbitem1}. The fact that equality holds if and only if $\vec{y} \in \Fix P_S \circ T$ is an immediate consequence of Proposition~\ref{prop:fixedpoints}.
\end{proof}

\subsection{Convergence}\label{sec:convergence}

In this subsection, we will show that sequences admitted by \GPPPPA \; are globally convergent to a point in $\Fix P_S \circ T$, whenever the latter is nonempty. The key result, Theorem~\ref{thm:PTQNE}, uses the following auxiliary lemma.

\begin{lemma}\label{lem:beta}
	Let $\Fix P_S \circ T \neq \emptyset$. The following holds: 
	\begin{equation*}%\label{Tycutter}
	(\forall \vec{x} \in \Fix P_S \circ T) (\forall \vec{y} \in S)\quad	\langle T\vec{x}-T\vec{y},\vec{y}-T\vec{y} \rangle \leq 0.
	\end{equation*}
\end{lemma}
\begin{proof}
	Let $\vec{x} \in \Fix P_S \circ T$ and $\vec{y} \in S$. By Lemma~\ref{lem:subdifferentialkappa}, we may consider two cases: when $T\vec{y} = P_{\cdk}\vec{y}$ and when $T\vec{y} = P_{\kappa \leq r_2}\vec{y}$ where $r_2 = \|T\vec{y}-\vec{y}\|$.
	
	\textbf{Case 1}: Let $T\vec{y} = P_{\cdk}\vec{y}$. Since $P_{\cdk}$ is FQNE, we have from Lemma~\ref{lem:BCQNE} that
	\begin{equation}\label{condition5}
	(\forall \vec{u} \in \Fix P_{\cdk})(\forall \vec{v} \in V)\quad \langle \vec{u}-P_{\cdk}\vec{v},\vec{v}-P_{\cdk}\vec{v} \rangle \leq 0.
	\end{equation}
	In particular $T\vec{x} \in \Fix P_{\cdk}$, and so we can apply \eqref{condition5} with $\vec{u} = T\vec{x}$ and $\vec{v} = \vec{y}$, obtaining
	\begin{equation}\label{condition5b}
	\langle T\vec{x}-P_{\cdk}\vec{y},\vec{y}-P_{\cdk}\vec{y} \rangle \leq 0.
	\end{equation}
	Finally, substituting in \eqref{condition5b} using the fact that $P_{\cdk}\vec{y} = T\vec{y}$, we obtain 
	\begin{equation*}
	\langle T\vec{x}-T\vec{y},\vec{y}-T\vec{y} \rangle \leq 0.
	\end{equation*}
	This shows the result in Case 1.
	
	\textbf{Case 2}: Let $T\vec{y} = P_{\kappa \leq r_2}\vec{y}$ where $r_2 = \|T\vec{y}-\vec{y}\|$. As $P_{\kappa \leq r_2}$ is FQNE, we have that 
	\begin{equation}\label{condition5c}
	(\forall \vec{u} \in \Fix P_{\kappa \leq r_2})(\forall \vec{v} \in V)\quad \langle \vec{u}-P_{\kappa \leq r_2}\vec{v},\vec{v}-P_{\kappa \leq r_2}\vec{v} \rangle \leq 0.
	\end{equation}
	We have by Lemma~\ref{lem:fgreaterthanb} that $r_2 \geq r' = \|T\vec{x}-\vec{x}\|$. Combining this with the fact from \cite[Theorem~4.4(ii)]{friedlander2019polar} that $\kappa (T\vec{x}) \leq \|T\vec{x}-\vec{x}\|$ , we have that $T\vec{x} \in \lev_{\leq r'} \kappa \subset \lev_{\leq r_2}\kappa$. Thus $T\vec{x} \in \Fix P_{\kappa \leq r_2}$, and so we can apply \eqref{condition5c} with $\vec{u}=T\vec{x}$ and $\vec{v}=\vec{y}$, obtaining
	\begin{equation}\label{condition5d}
	\langle T\vec{x}-P_{\kappa \leq r_2}\vec{y},\vec{y}-P_{\kappa \leq r_2}\vec{y} \rangle \leq 0.
	\end{equation}
	Finally, since $P_{\kappa \leq r_2}\vec{y} = T\vec{y}$, we may substitute in \eqref{condition5d} to obtain
	\begin{equation*}
	\langle T\vec{x}-T\vec{y}, \vec{y}-T\vec{y} \rangle \leq 0.
	\end{equation*}
	This concludes the result.	
\end{proof}

\begin{fact}\label{fact:subspacepythagorean}
	Since $S$ is an affine subspace, it holds that
	\begin{equation}\label{subspacepythagorean}
	(\forall \vec{u} \in X\times \R)(\forall \vec{v} \in S)\quad \|\vec{u}-\vec{v}\|^2 = \|P_S(\vec{u})-\vec{v}\|^2+\|\vec{u}-P_S(\vec{u})\|^2.
	\end{equation}
\end{fact}
\begin{proof}
	This follows immediately from the Pythagorean theorem.
\end{proof}

\begin{figure}
	\begin{center}
		\begin{tikzpicture}[scale=7.0]
		%draw the set S		
		\draw [gray, thin,dashed] (1,0) -- (1,-0.85); 
		\node [above,gray] at (1,0) {$S=X \times \{1\}$};
		
		\draw [fill,black] (1,-0.1) circle [radius=0.01];
		\node [right] at (1,-0.1) {$\vec{x}=P_S \circ T\vec{x}$};
		
		\draw [fill,gray] (0,0) circle [radius=0.01];
		\node [gray, left] at (0,0) {$\vec{0}$};
		
		%inner level curve
		\draw [red, dotted] (0.1,0) -- (0.1,-0.2);
		\draw [red, dotted] (0.1,-0.2) -- (-0.2,-0.2-.25*.3); %slope 1/4
		
		%outer level curve
		\draw [red,dotted] (0.2,0) -- (0.2,-0.4);
		\draw [red,dotted] (0.2,-0.4) -- (-0.2,-0.4-0.1);  %slope 1/4
		
		\draw [fill,black] (0.1,-0.1) circle [radius=0.01];
		\node [above right] at (0.1,-0.1) {$T\vec{x}$};
		
		%connect x and Tx
		\draw [<->] (0.15,-0.1) -- (0.95,-0.1);
		
		%y, Ty, and PSoTy
		\draw [fill,black] (1,-0.8) circle [radius=0.01];
		\node [right] at (1,-0.8) {$\vec{y}$};
		\draw [fill,black] (0.2,-0.4) circle [radius=0.01];
		\node [below] at (0.2,-0.4) {$T\vec{y}$};
		\draw [fill,black] (1,-0.4) circle [radius=0.01];
		\node [right] at (1,-0.4) {$P_S \circ T\vec{y}$};
		\draw [->] (0.25,-0.4) -- (0.95,-0.4);
		
		%connect y and Ty
		\draw [->] (0.95,-0.775) -- (0.25,-0.425);
		
		%convex separation for Ty and Tx
		\draw [gray, dashed] (0,-0.8) -- (0.4,0);
		
		%draw my added point Tx-u
		\draw [fill,blue] (0,-0.3) circle [radius=0.01];
		\node [left,blue] at (0,-0.3) {$T\vec{x}-\vec{u}$};
		
		%draw my extra triangle
		\draw [blue,->] (0.025,-0.3+.05) -- (0.1-0.025,-0.1-0.05);
		\draw [blue] (.1+0.025, -0.1-.05) -- (0.2-0.025,-0.4+0.05);
		\draw [blue,<-] (0.05,-0.3-0.025) -- (0.2-0.05,-0.4 + 0.025);
		
		%draw my extra hypotenuse
		\draw [blue] (0.1+.05,-0.1-0.04) -- (0.95,-0.8+0.04);
		
		%label my a and c segments and label them
		\draw [purple,|-|] (1,-0.15) -- (1,-0.75);
		\node [right,purple] at (1,-0.6) {$\|\vec{y}-\vec{x}\|$};
		
		\draw [|-|,blue] (1.05,-0.15) -- (1.05,-.35);
		\node [right,blue] at (1.05,-0.25) {$\|\vec{x}-P_S\circ T\vec{y} \|$};
		
		%label my extra segments
		\node [above] at (0.5,-0.1) {$\|T\vec{x}-\vec{x}\|$};
		\node [above right,blue] at (.69,-0.56) {$\|T\vec{x}-\vec{y}\|$};
		\node [below left] at (0.6,-0.6) {$\|T\vec{y}-\vec{y}\|$};
		\node [above left,blue] at (0.05,-0.2) {$\vec{u}$};
		\node [below left,blue] at (0.1,-0.35) {$\beta(T\vec{y}-\vec{y})$};
		%\node [right,blue]  at (0.15,-0.25) {$Tx-Ty$};

		\end{tikzpicture}
	\end{center}
	\caption{The strategy for Theorem~\ref{thm:PTQNE} is illustrated.}\label{fig:PTQNE}
\end{figure}
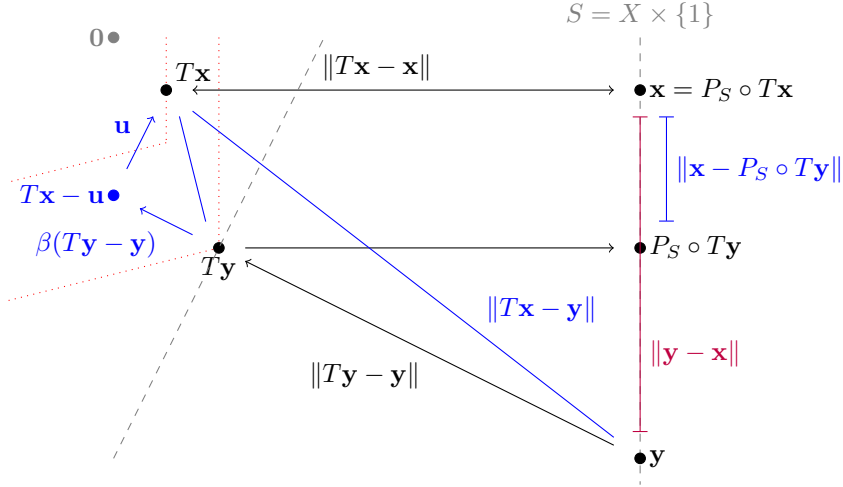

Figure~\ref{fig:PTQNE} illustrates the strategy of the following theorem, which brings together all of the different results we have established so far.

\begin{theorem}\label{thm:PTQNE}
	Let $\Fix P_S \circ T \neq \emptyset$. Set $\vec{y}_0 \in S$. The operator $P_S \circ T: S \rightarrow S$ is SQNE, and so the sequence $(\vec{y}_n)_{n \in \N} \subset S$ given by 
	\begin{equation*}%\label{eqn:Fejer}
		\vec{y}_{n+1}:=P_S \circ T \vec{y}_n
	\end{equation*}
	is Fej\'er monotone with respect to $\Fix (P_{S} \circ T)$. More specifically,
	\begin{equation*}
	(\forall \vec{x} \in \Fix P_S \circ T)(\forall \vec{y} \in S)\quad \|P_S \circ T\vec{y} -\vec{x}\|^2 \leq \|\vec{y}-\vec{x}\|^2 -E(\vec{x},\vec{y}),%\label{criticalinequality}
	\end{equation*}
	where $E$ is as in \eqref{def:E}.
\end{theorem}
\begin{proof}
	Let $\vec{x} \in \Fix P_S \circ T$ and $\vec{y} \in S$. There exists $\beta \in \R$ such that 
	\begin{equation}\label{constructingtriangle}
		T\vec{x}-\vec{y} = (1+\beta)(T\vec{y}-\vec{y})+\vec{u}\quad\text{for some}\; \vec{u} \in ({\rm span} \{T\vec{y}-\vec{y}\})^\perp.
	\end{equation}
	We will first show that $\beta \geq 0$. From Lemma~\ref{lem:beta} we have that 
	\begin{equation}\label{condition5beta1}
		\langle T\vec{x}-T\vec{y},\vec{y}-T\vec{y} \rangle \leq 0.
	\end{equation}
	Adding $\langle \vec{y}-T\vec{x},\vec{y}-T\vec{y} \rangle$ to both sides of \eqref{condition5beta1} yields	
	\begin{equation}\label{condition5beta}
		\|T\vec{y}-\vec{y}\|^2 \leq \langle T\vec{x}-\vec{y},T\vec{y}-\vec{y} \rangle.
	\end{equation}
	Combining \eqref{constructingtriangle} and \eqref{condition5beta}, we obtain
	\begin{align}
		\|T\vec{y}-\vec{y}\|^2 &\leq \langle (1+\beta)(T\vec{y}-\vec{y})+\vec{u},T\vec{y}-\vec{y} \rangle\nonumber \\
		&=(1+\beta)\|T\vec{y}-\vec{y}\|^2.\label{whybetageqzero}
	\end{align}
	Now \eqref{whybetageqzero} implies that $\beta \geq 0$ or $\|T\vec{y}-\vec{y}\| = 0$. If $\|T\vec{y}-\vec{y}\| = 0$, then $T\vec{y}=\vec{y} \in S$ and so $P_S\circ T \vec{y} = \vec{y}$, and so $\vec{y} \in \Fix (P_S \circ T)$, in which case we are done. Thus we may restrict to considering the case when $\|T\vec{y}-\vec{y}\| > 0$. This, together with \eqref{whybetageqzero}, yields
	\begin{equation*}%\label{betageqzero}
		\beta \geq 0. 
	\end{equation*}
	Now applying the Pythagorean Theorem to \eqref{constructingtriangle} yields
	\begin{equation}\label{Pythag1}
		\|\vec{u}\|^2 = \|T\vec{x}-\vec{y}\|^2- \|(1+\beta)(T\vec{y}-\vec{y})\|^2.
	\end{equation}
	Moreover, we may rearrange \eqref{constructingtriangle} to obtain
	\begin{equation}\label{Pythag2a}
		T\vec{x}-T\vec{y} = \beta(T\vec{y}-\vec{y})+\vec{u}. 
	\end{equation}
	Applying the Pythagorean Theorem to \eqref{Pythag2a}, we obtain 
	\begin{equation}\label{Pythag2}
		\|T\vec{x}-T\vec{y}\|^2 = \|\beta(T\vec{y}-\vec{y})\|^2 + \|\vec{u}\|^2.
	\end{equation}
	Using \eqref{Pythag1} to substitute for $\|\vec{u}\|^2$ in \eqref{Pythag2}, we obtain
	\begin{align}
		\|T\vec{x}-T\vec{y}\|^2 &= \|\beta(T\vec{y}-\vec{y})\|^2 + \|T\vec{x}-\vec{y}\|^2- \|(1+\beta)(T\vec{y}-\vec{y})\|^2 \nonumber\\
		&= \|T\vec{x}-\vec{y}\|^2-(1+2\beta)\|T\vec{y}-\vec{y}\|^2. \label{BigA}
	\end{align}	
	Now from Fact~\ref{fact:subspacepythagorean} we have that \eqref{subspacepythagorean} holds for $\vec{u}=T\vec{x}$ and $v=\vec{y}$, and so we have
	\begin{align}
		\|T\vec{x}-\vec{y}\|^2 &= \|T\vec{x}-P_S\circ T(\vec{x})\|^2 + \|\vec{y}-P_S\circ T(\vec{x})\|^2\nonumber \\
		\text{which is}\quad \|T\vec{x}-\vec{y}\|^2 &= \|T\vec{x}-\vec{x}\|^2 + \|\vec{y}-\vec{x}\|^2.\label{Pythag3}
	\end{align}
	Now we may use \eqref{Pythag3} to substitute for $\|T\vec{x}-\vec{y}\|^2$ in \eqref{BigA} and obtain
	\begin{equation}\label{BigB}
		\|T\vec{x}-T\vec{y}\|^2 = \|T\vec{x}-\vec{x}\|^2 + \|\vec{y}-\vec{x}\|^2 -(1+2\beta)\|T\vec{y}-\vec{y}\|^2.
	\end{equation}
	Now we have that
	\begin{equation}\label{E1}
	\|T\vec{y}-\vec{y}\|^2 =  \|T\vec{x}-\vec{x}\|^2 + \|T\vec{y}-\vec{y}\|^2-\|T\vec{x}-\vec{x}\|^2 = \|T\vec{x}-\vec{x}\|^2 + E(\vec{x},\vec{y}),
	\end{equation}
	where $E$ is as defined in \ref{def:E}. Multiplying both sides of \eqref{E1} by $-(1+2\beta)$, we obtain
	\begin{equation}\label{E2}
	-(1+2\beta)\|T\vec{y}-\vec{y}\|^2 = -(1+2\beta)\|T\vec{x}-\vec{x}\|^2 -(1+2\beta)E(\vec{x},\vec{y}).
	\end{equation}
	Using \eqref{E2} to make the appropriate substitution for $-(1+2\beta)\|T\vec{y}-\vec{y}\|^2$ in \eqref{BigB}, we obtain
	\begin{align}
	\|T\vec{x}-T\vec{y}\|^2 &= \|T\vec{x}-\vec{x}\|^2 + \|\vec{y}-\vec{x}\|^2 -(1+2\beta)\|T\vec{x}-\vec{x}\|^2 -(1+2\beta)E(\vec{x},\vec{y})\nonumber \\
	&= \|\vec{y}-\vec{x}\|^2 -2\beta\|T\vec{x}-\vec{x}\|^2 -(1+2\beta)E(\vec{x},\vec{y}).\label{E3}
	\end{align}
	Since $S$ is closed and convex, $P_S$ is nonexpansive (\cite[Proposition~4.16]{BC}). Thus we have that
	\begin{equation}\label{aleqv}
	\|P_S \circ T\vec{y} - P_S \circ T\vec{x}\|^2 \leq \| T\vec{y}-T\vec{x}\|^2.
	\end{equation}
	Since $\vec{x} \in \Fix P_S \circ T$, we have that $P_S \circ T \vec{x} = \vec{x}$. Making this substitution in \eqref{aleqv}, we obtain
	\begin{equation}\label{E4}
	\|P_S \circ T\vec{y} -\vec{x} \|^2 \leq \|T\vec{y}-T\vec{x}\|^2.
	\end{equation}
	Together, \eqref{E3} and \eqref{E4} yield
	\begin{align}
	\|P_S \circ T\vec{y} -\vec{x} \|^2 &\leq \|\vec{y}-\vec{x}\|^2 -2\beta\|T\vec{x}-\vec{x}\|^2 -(1+2\beta)E(\vec{x},\vec{y})\\
	&\leq \|\vec{y}-\vec{x}\|^2-E(\vec{x},\vec{y}).\label{E5}
	\end{align}
	where the second inequality uses the fact that $\beta \geq 0$ and $E(\vec{x},\vec{y}) \geq 0$. This shows the desired result.
\end{proof}

Theorem~\ref{thm:PTQNE} admits the following corollary.

\begin{corollary}[Averaged variant]\label{cor:averagedFQNE}
	Let $\Fix P_S \circ T \neq \emptyset$. The operator given by
	\begin{equation*}
	\frac12 P_S \circ T + \frac12 \Id
	\end{equation*}
	is FQNE. 
\end{corollary}
\begin{proof}
	By Theorem~\eqref{thm:PTQNE}, we have that $P_S \circ T$ is QNE. By Lemma~\ref{lem:BCQNE}, an operator $R$ is FQNE if and only if $2U-\Id$ is QNE. Letting
	\begin{equation*}
	2U-\Id = P_S \circ T,
	\end{equation*}
	we have that $U$ is FQNE and that $U=\frac12 P_S \circ T + \frac12 \Id$.
\end{proof}

Now having the key results of Theorem~\ref{thm:PTQNE}, we are ready to show convergence for both \GPPPPA \;and its under-relaxed variants.

\begin{theorem}[Convergence of under-relaxed variants of \GPPPPA]\label{thm:PTrelaxedstrong}
	Let $\Fix P_S \circ T \neq \emptyset$. Let $\gamma \in \left]0,1 \right[$ and $\vec{y}_0 \in S$. The sequence given by 
	\begin{align*}
	\vec{y}_{n+1}:=&\;\mathcal{U}_\gamma \vec{y}_n,\\
	\text{where}\quad \mathcal{U}_\gamma:=&\;(1-\gamma)P_S \circ T +\gamma \Id
	\end{align*}
	is strongly convergent to some $\vec{y} \in \Fix P_S \circ T$. 
\end{theorem}
\begin{proof}
	Notice that 
	\begin{align*}
	\mathcal{U}_\gamma&=((1-\gamma)(P_S \circ T -\Id)+\Id) \\ 
	&=\left(1-\frac{2\gamma}{2}\right)(P_S \circ T-\Id)+\Id\\
	&=(2-2\gamma)\left(\frac{1}{2}P_S\circ T - \frac12 \Id \right)+\Id\\
	&=(2-2\gamma)\left(\left(\frac12 P_S\circ T+\frac12 \Id \right)-\Id \right) + \Id\\
	&=(2-2\gamma)(U-\Id)+\Id,
	\end{align*}
	where $2\gamma \in \left]0,2\right[$ and 
	\begin{equation*}
	U:=\frac12 P_S \circ T + \frac12 \Id
	\end{equation*}
	is the FQNE operator from  Corollary~\ref{cor:averagedFQNE}. Thus, applying Lemma~\ref{lem:DLRcutterresult} for the operator $U$, we have that $\mathcal{U}_\gamma$ is QNE and that 
	\begin{align}
	\|P_S \circ T\vec{y}_n - \vec{y}_n \|&= 2\left \|\frac12 P_S \circ T\vec{y}_n - \frac12 \vec{y}_n \right \|\nonumber \\
	&= 2\left \|\vec{y}_n-\left(\frac12 P_S \circ T + \frac12 \Id \right)\vec{y}_{n}\right \|\nonumber \\
	&=2\|\vec{y}_n-U(\vec{y}_n)\| \rightarrow 0 \quad \text{as}\quad n \rightarrow \infty.\label{triangleinequality1}
	\end{align}
	Since $(\vec{y}_n)_{n \in \N}$ is Fej\'er monotone, it is bounded. Since $S$ is finite dimensional and $(\vec{y}_n)_{n \in \N} \subset S$ is bounded, we may take a convergent subsequence $(\vec{y}_j)_{j \in J \subset \N}$ such that
	\begin{equation}\label{triangleinequality2}
	\|\vec{y}_{j} - \vec{y}\| \rightarrow 0 \quad \text{as}\; j \rightarrow \infty
	\end{equation}
	for some $\vec{y} \in S$. Since $T$ is continuous and $P_S$ is continuous, $P_S\circ T$ is continuous. Thus we have
	\begin{equation}\label{triangleinequality0}
	\underset{j \rightarrow \infty}{\lim}P_S \circ T(\vec{y}_{j}) = P_S \circ T(\underset{j \rightarrow \infty}{\lim}\vec{y}_{j}) = P_S \circ T (\vec{y}).
	\end{equation}
	%Finally, notice that by \eqref{triangleinequality0},
	%\begin{equation}\label{triangleinequality3}
	%\| y_{j}-y_{j+1} \| = (1-\gamma)\|y_j-T\circ P_S\circ T(y_j) \|  \rightarrow 0 \quad \text{as}\; n_k \rightarrow \infty.
	%\end{equation}
	The triangle inequality yields
	\begin{align}\label{triangleinequality4}
	\|\vec{y}-P_S\circ T(\vec{y})\| \leq \|\vec{y}-\vec{y}_{j}\|+ \|\vec{y}_{j}-P_S\circ T(\vec{y}_{j})\|+\|P_S\circ T(\vec{y}_{j})-P_S\circ T(\vec{y})\|.
	\end{align}
	Taking the limit as $j \rightarrow \infty$, each of the terms in the right hand side of \eqref{triangleinequality4} go to zero by \eqref{triangleinequality2}, \eqref{triangleinequality1}, and \eqref{triangleinequality0} respectively. Thus $\|\vec{y}-P_S \circ T(\vec{y})\| = 0$, and so $\vec{y} = P \circ T(\vec{y})$. Thus we have that $\vec{y} \in \Fix P_S \circ T$. 
	
	Since $(\vec{y}_n)_{n \in \N}$ is Fej\'er monotone with respect to $\Fix P_S \circ T$ and possesses a sequential cluster point $\vec{y} \in \Fix P_S \circ T$, we conclude by Theorem~\ref{thm:baucomstrongconvergence} that $\vec{y}_n \rightarrow \vec{y}$ as $n \rightarrow \infty$. 
\end{proof}

Having proven the convergence for the under-relaxed variants of \GPPPPA, we now show the convergence of its non-relaxed version.

\begin{theorem}[Convergence of \GPPPPA]\label{thm:PTstrong}
	Let $\Fix P_S \circ T \neq \emptyset$. Let $\vec{y}_0 \in S$. The sequence $(\vec{y}_n)_{n \in \N}$ given by 
	\begin{equation*}
	\vec{y}_{n+1}:= P_S \circ T \vec{y}_n
	\end{equation*}
	converges to a point $\vec{y} \in \Fix P_S \circ T$. 
\end{theorem}
\begin{proof}
	Fix $\vec{x} \in \Fix P_S \circ T$. Applying Theorem~\ref{thm:PTQNE}, we have that 
	\begin{equation*}
	\|\vec{y}_{n+1}-\vec{x}\|^2 \leq \|\vec{y}_n-\vec{x}\|^2 - E(\vec{x},\vec{y}_n),
	\end{equation*}
	and so we have that
	\begin{equation*}
	0 \leq \|\vec{y}_{n+1}-\vec{x}\|^2 \leq \|\vec{y}_0-\vec{x}\|^2 -\sum_{i=0}^n E(\vec{x},\vec{y}_i).
	\end{equation*}
	Thus we obtain
	\begin{equation*}
	\sum_{i=0}^n E(\vec{x},\vec{y}_i) \leq \|\vec{y}_0-\vec{x}\|^2,
	\end{equation*}
	which shows that
	\begin{equation}\label{Etozero}
	E(\vec{x},\vec{y}_n) \rightarrow 0 \quad \text{as}\; n \rightarrow \infty.
	\end{equation}
	From the definition of $E$, \eqref{Etozero} implies that
	\begin{equation}\label{downtor1}
	\|T\vec{y}_n - \vec{y}_n\| \downarrow \|T\vec{x}-\vec{x}\| \quad\text{as}\; n\rightarrow \infty.
	\end{equation}
	As $(\vec{y}_n)_{n \in \N}$ is Fej\'er monotone, it is bounded. Thus we may take a convergent subsequence $(\vec{y}_j)_{j \in J \subset \N}$. Therefore, let
	\begin{equation*}
	\vec{y}_j \rightarrow \vec{y} \quad \text{as}\quad j \rightarrow \infty.
	\end{equation*}
	Combining with \eqref{downtor1}, we obtain
	%\begin{equation}
	%\|Ty-y\|  = \left \|\underset{j \rightarrow \infty}{\lim}(T-\Id)y_j  \right \| = \underset{j \rightarrow \infty}{\lim}\|Ty_j - y_j \| = \|Tx-x\|,
	%\end{equation}
	\begin{equation*}
	\|T\vec{y}-\vec{y}\|   = \underset{j \rightarrow \infty}{\lim}\|T\vec{y}_j - \vec{y}_j \| = \|T\vec{x}-\vec{x}\|,
	\end{equation*}
	where the first equality follows from the continuity of $T-\Id$ and the second equality is from \eqref{downtor1}. Now since $\|T\vec{y}-\vec{y}\| = \|T\vec{x}-\vec{x}\|$ with $\vec{x} \in \Fix P_S \circ T$, we have by Proposition~\ref{prop:fixedpoints} that $\vec{y} \in \Fix P_S \circ T$. Since $(\vec{y}_n)_{n \in \N}$ is Fej\'er monotone with respect to $\Fix P_S \circ T$ and possesses a sequential cluster point $\vec{y} \in \Fix P_S \circ T$, we conclude by Theorem~\ref{thm:baucomstrongconvergence} that $\vec{y}_n \rightarrow \vec{y}$ as $n \rightarrow \infty$. 
\end{proof}

When $\Fix P_S \circ T \neq \emptyset$, Theorem~\ref{thm:PTstrong} guarantees the convergence of \GPPPPA, and Theorem~\ref{thm:PTrelaxedstrong} does the same for its under-relaxed variants. The next corollary simply formalizes this by including both cases.

\begin{corollary}\label{cor:relaxedandunrelaxed}
	Let $\Fix P_S \circ T \neq \emptyset$. Let $\gamma \in \left[0,1 \right[$ and $\vec{y}_0 \in S$. Then the sequence given by 
	\begin{align}
	\vec{y}_{n+1}:=&\mathcal{U}_\gamma \vec{y}_n \nonumber\\
	\text{where}\quad \mathcal{U}_\gamma:=&(1-\gamma)P_S \circ T +\gamma \Id \label{relaxedorunrelaxed}
	\end{align}
	is convergent to some $\vec{y} \in \Fix P_S \circ T$, and the operator $\mathcal{U}_\gamma$ is paracontracting.
\end{corollary}
\begin{proof}
	The convergence when $\gamma \in \left]0,1\right[$ is shown by Theorem~\ref{thm:PTrelaxedstrong}, and the convergence when $\gamma = 0$ is dealt with by Theorem~\ref{thm:PTstrong}. These theorems also show that $\mathcal{U}_\gamma$ is QNE in these two cases respectively. Combining with the fact that $\mathcal{U}_\gamma$ is obviously a weighted average of continuous operators, the paracontracting property is clear.
\end{proof}

\subsection{Shadow sequence behaviour}\label{sec:shadow}

Having established convergence for the governing sequence, we also have the following result that describes the behaviour of the sequence of shadows of the proximity operator: $T\circ(\mathcal{U}_\gamma)^{n-1}(y_0,1)$.

\begin{corollary}\label{cor:lambdashadow}
	Let $\Fix P_S \circ T \neq \emptyset$. Set $(y_0,1) \in S$ and $\gamma \in \left[0,1\right[$. The sequence $(\lambda_n)_{n \in \N} \subset S$ given by 
	\begin{equation*}%\label{lambdaconvergence}
	(y_{n+1},\lambda_{n+1}):=T\circ(\mathcal{U}_\gamma)^{n}(y_0,1),
	\end{equation*}
	where $\mathcal{U}_\gamma$ is as specified in \eqref{relaxedorunrelaxed} satisfies 
	\begin{equation*}
	\lambda_n \rightarrow 1-r' \in \left[0,1\right] \quad \text{as}\; n \rightarrow \infty,
	\end{equation*}
	where $r'$ is as specified in \eqref{characterizationofr1}. 
\end{corollary}
\begin{proof}
	From Corollary~\ref{cor:relaxedandunrelaxed} we have that
	$$
	\vec{x}_n:=(y_n,1) \rightarrow \vec{x} \;\;\text{for some}\;\; \vec{x} \in \Fix P_S \circ T.
	$$
	Given that $\vec{x} =(y,1)$ for some $y \in X$, we clearly have
	\begin{equation*}
		(y_n,1) \rightarrow (y,1).
	\end{equation*}
	From Proposition~\ref{prop:fixedpoints}, we have that
	\begin{equation}\label{lambda1}
	\|T(y,1)-(y,1)\| = r',
	\end{equation}
	where $r'$ is as characterized in \eqref{characterizationofr1}. Using Lemma~\ref{lem:abslambdaleq1}, we have that $T(y,1)=(y,\lambda)$ for some $\lambda \in \left[0,1\right]$. This yields
	\begin{equation}\label{lambda2}
	\|T(y,1)-(y,1)\| = \|(y,\lambda)-(y,1)\| = 1-\lambda.
	\end{equation}
	Combining \eqref{lambda1} and \eqref{lambda2} we have that $r' = 1-\lambda$, and so
	\begin{equation*}
	\lambda = 1-r' \in \left[0,1\right].
	\end{equation*}
	Since $(y_n,1) \rightarrow (y,1)$ and $T$ is continuous, we have that 
	$$
	(y_{n+1},\lambda_{n+1}) = T(y_n,1) \rightarrow T(y,1) = (y,\lambda).
	$$ 
	Thus $\lambda_n \rightarrow \lambda = 1-r'$. This concludes the result.
\end{proof}

\subsection{The operator $P_S \circ T$ is not, generically, FQNE}\label{sec:counterexample}

The property of firm quasinonexpansivity is especially important in the analysis of algorithms. In this section, we discuss under which conditions the operator $P_S \circ T$ may or may not exhibit this property. In particular, we provide an example illustrating that it is not, generically, FQNE. First we show, in Proposition~\ref{prop:PTalmostcutter}, that failure to be a FQNE operator implies some specific conditions.

\begin{proposition}\label{prop:PTalmostcutter}
	Let $(y,1) \in S$, and let $(x,1) \in \Fix P_S \circ T$. Let $(x,\lambda) = T(x,1)$ and $(w,\mu)=T(y,1)$. If
	\begin{equation}\label{notacutter}
	0 < \langle (x,1)-P_S \circ T(y,1),(y,1)-P_S \circ T(y,1) \rangle,
	\end{equation}
	then the following hold:
	\begin{enumerate}[label=(\roman*)]
		\item\label{lambdamuinequality} $\lambda < \mu < 1$;
		\item\label{TxisnotPdomx} $T(x,1) \neq P_{\cdk}(x,1)$.
	\end{enumerate}
\end{proposition}
\begin{proof}
	Suppose that \eqref{notacutter} holds. Then we have that
	\begin{align}
	0 &< \langle (x,1)-P_S \circ T(y,1),(y,1)-P_S \circ T(y,1) \rangle \label{PTC3}\\
	 &=\langle (x,1) - (w,1), (y,1)-(w,1) \rangle \nonumber \\
	&= \langle (x-w,0),(y-w,0) \rangle  \nonumber \\
	&= \langle x-w,y-w \rangle. \label{PTC0}
	\end{align}
	
	We will first show \ref{lambdamuinequality}.\\
	
	From Lemma~\ref{lem:beta} we have that
	\begin{equation}\label{PTC1}
	\langle T(x,1)-T(y,1), (y,1)-T(y,1) \rangle \leq 0.
	\end{equation}
	Let $(w,\mu) = T(y,1)$ and $(x,\lambda)=T(x,1)$ for $\lambda \leq 1$. Then \eqref{PTC1} becomes
	\begin{align}
	\langle (x,\lambda)-(w,\mu),(y,1)-(w,\mu) \rangle &\leq 0\nonumber \\
	\langle (x-w,\lambda-\mu), (y-w,1-\mu) \rangle &\leq 0\nonumber  \\
	\langle x-w,y-w \rangle + (\lambda-\mu)(1-\mu) &\leq 0. \label{PTC2}
	\end{align}
	By Lemma~\ref{lem:abslambdaleq1}\ref{eqn:abslambdaleq1}, $\lambda \in \left[0,1\right]$ and so there are only four possibilities: when $\lambda =1$, when $\mu \leq \lambda <1$, when $\lambda < 1 \leq \mu$, and when $\lambda < \mu < 1$. We will show that any case other than $\lambda < \mu < 1$ implies a contradiction. 
	
	\textbf{Case $\lambda =1$}. Suppose $\lambda =1$. Then $(\lambda-\mu)(1-\mu) = (1-\mu)^2 \geq 0$. 
	
	Combining this fact with \eqref{PTC2}, we obtain
	\begin{align}\label{contradictPTC0}
	\langle x-w,y-w \rangle \leq 0,
	\end{align}
	which contradicts \eqref{PTC0}, and so we obtain a contradiction. This concludes the case $\lambda=1$.
	
	\textbf{Case $\mu \leq \lambda < 1$}. Suppose $\mu \leq \lambda \leq 1$. We have that $(\lambda-\mu) \geq 0$ and $(1-\mu)\geq 0$, and so clearly $(\lambda-\mu)(1-\mu) \geq 0$. Combining this with \eqref{PTC2}, we again obtain \eqref{contradictPTC0}, which is a contradiction. This concludes the case when $\mu \leq \lambda < 1$.
	
	\textbf{Case $\lambda < 1 \leq \mu$}. Suppose $\lambda < 1 \leq \mu$. Then $(\lambda -\mu) \leq 0$ and $(1-\mu) \leq 0$, and so $(\lambda-\mu)(1-\mu) \geq 0$. Combining this with \eqref{PTC2}, we again obtain \eqref{contradictPTC0}, a contradiction. This concludes the case when $\lambda < 1 \leq \mu$.
	
	We are left with only one possibility, $\lambda < \mu < 1$, and so \ref{lambdamuinequality} holds. \\
	
	We next show \ref{TxisnotPdomx}. Having established that $\lambda < \mu$, we have by Lemma~\ref{lem:fgreaterthanb}\ref{lem:fgreaterthanbitem2c} that $(\lambda < \mu) \implies T(x,1) \neq P_{\cdk}(x,1)$. Thus we have \ref{TxisnotPdomx}.
\end{proof}
Proposition~\ref{prop:PTalmostcutter} shows that any example of a gauge for which $P_S \circ T$ fails to be firmly quasi-nonexpansive must satisfy both \ref{lambdamuinequality} and \ref{TxisnotPdomx}. Now we will see an example that satisfies both of these properties and serves as a counter-example to the tempting idea that $P_S \circ T$ is generically FQNE.
\begin{figure}[ht]
	\begin{center}
		\includegraphics[width=.49\textwidth]{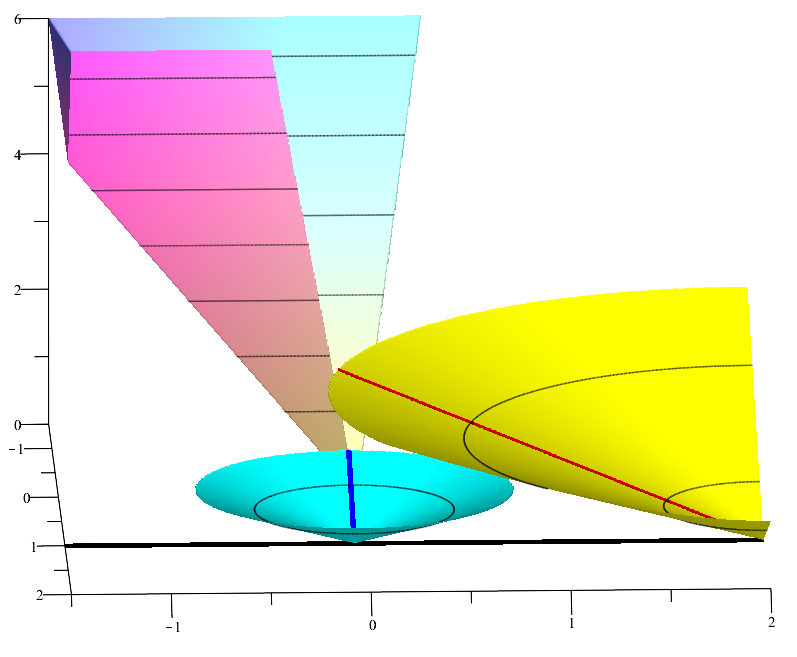}\;\includegraphics[width=.49\textwidth]{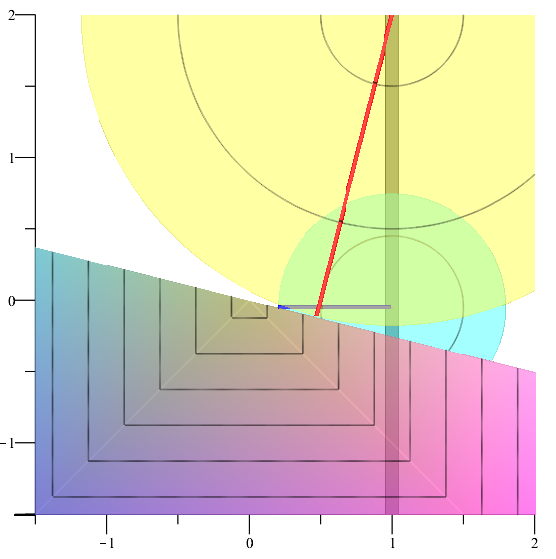}
	\end{center}
	\caption{An illustration of Example~\ref{ex:cutter_counterexample}.}\label{fig:cutter_counterexample}
\end{figure}
	
\begin{example}[$P_S \circ T$ is generically not a cutter]\label{ex:cutter_counterexample}
	Let 
	$$\kappa : \R^1 \times \R \rightarrow \R: \vec{x} \mapsto 4\| \vec{x} \|_\infty + \iota_C \vec{x},$$
	where $C:= \{(v,u)\;|\; v\leq-u/4 \}$. Then 
	$$
	T_\gamma(1,-1/20) = P_{\kappa \leq 1/5}(1,-1/20) = (1/5,-1/20),
	$$ 
	and so $(1,-1/20) \in \Fix P_S \circ T_\gamma$. Additionally, 
	$$
	T_\gamma(2,1) = P_{\cdk}(2,1) = (-2/17,8/17),
	$$
	and so $(2,1) \notin \Fix P_S \circ T_{\gamma}$. We have that
	\begin{equation*}
	\langle P_S \circ T_{\gamma}(2,1) - (-1/20,1) , (-1/20,1) - (2,1) \rangle = (-2/17+1/20)(-1/20-2) > 0.
	\end{equation*}
	This example is illustrated in Figure~\ref{fig:cutter_counterexample}.
\end{example}

\section{Fundamental set and existence of fixed points}\label{sec:Fixnonempty}

In the previous section, we consistently assumed that $\Fix P_S \circ T \neq \emptyset$. It bears noting that this condition may not hold for a \textit{general gauge}. Of course, it \textit{does} hold for \PPPPA, under the conditions in Theorem~\ref{thm:fixedpoints_argmin}. We will provide sufficient conditions for the more general \GPPPPA, which allow us to describe the solutions to \PPPPA \; as lying on an exposed face of a dilated fundamental set. For the purpose, we make use of the Minkowski function representation of the gauge:
\begin{equation}\label{eqn:Minkowski}
\kappa = \gamma_D :X \times \R \rightarrow \R: \vec{x}\mapsto \inf \{\mu \geq 0 \;|\; \vec{x} \in \mu D \}.
\end{equation}
Such a representation always holds by choosing $D = \{\vec{x} \in X\times \R \;|\; \kappa(\vec{x}) \leq 1\}$ \cite{friedlander2019polar}.

The following Lemma will be instrumental to our main result in Theorem~\ref{thm:existence}.

\begin{lemma}\label{lem:conedom}
	Let $D= \lev_{\kappa \leq 1}$. The following hold.
	\begin{enumerate}[label=(\roman*)]
		\item\label{lem:conedom1} $\cone D = \dom \kappa$.
		\item\label{lem:conedom2} If there exists $\lambda' > 0$ such that $$
		\lambda' = \max_{\lambda \in \R}\{\lambda \;\;|\;\; \exists y' \in X \text{so\;that}\;(y',\lambda') \in D  \},$$
		then for any $(y,\beta) \in \cone D \cap (X\times \R_{>0})$ there exists a minimal $0<r<\infty$ so that $(y,\beta) = r \vec{d}$ for some $\vec{d} \in D$.
	\end{enumerate}
\end{lemma}
\begin{proof}
	\ref{lem:conedom1}: Let $\vec{x} \in \cone D$. Then there exist $\lambda <\infty, \vec{d} \in D$ such that $\vec{x} = \lambda \vec{d}$. By positive homogeneity, $\kappa(\lambda \vec{d})= \lambda \kappa(\vec{d}) \leq \lambda < \infty$, and so $\vec{x} \in \dom \kappa$, so $\cone D \subset \dom \kappa$. Now let $\vec{x} \in \dom \kappa$. Then there exists $r<\infty$ such that $\kappa(\vec{x}) = r$ and so by homogeneity $\kappa(\vec{x}/r) = 1$ and so $\vec{x}/r \in D$, so $\vec{x}= r(\vec{x}/r) \in \cone D$, and so $\dom \kappa \subset \cone D$. 
	
	\ref{lem:conedom2}: Let $(y,\beta) \in \cone D \cap (X\times \R_{>0})$. First of all, notice that by \ref{lem:conedom1}, $(y,\beta) \in \dom \kappa$ and so 
	$$
	r:=\inf\{\lambda \;|\; (y,\beta) \in \lambda D \}< \infty.
	$$
	Next we show $r>0$. Since $(y,\beta) \in \cone D$, there exists \emph{some} $\lambda\geq 0, (d_y,\mu)\in D$ such that $\lambda (d_y,\mu) = (y,\beta)$. Now any such $\lambda$ that satisfies this equality clearly satisfies $\lambda \mu = \beta$ with $\beta >0$, and so all three constants are greater than zero. Moreover, any such constant $\lambda$ that satisfies this equality satisfies $\lambda = \beta/\mu \geq \beta / \lambda' > 0$. Thus we have that
	$$
	r=\inf\{\lambda \;|\; (y,\beta) \in \lambda D \} \geq \beta/\lambda' > 0.
	$$
	Now let $(\lambda_n)_n$ satisfy $\lambda_n \downarrow r$ as $n \rightarrow \infty$. Since $\lambda_n >r$ for all $n$, $(y,\beta) \in \lambda_n D$ for all $n$ and so there exists $(\vec{d}_n)_n$ such that $(y,\beta)=\lambda_n \vec{d}_n$ for all $n$. Notice that
	$$
	\|\vec{d}_n\| = \frac{\|(y,\beta)\|}{\lambda_n} \leq \frac{\|(y,\beta)\|}{r}< \infty,
	$$
	and so the sequence $(\vec{d}_n)_n$ is bounded. Thus we can pass to a convergent subsequence if need be and have $\vec{d}_n \rightarrow \vec{d}$ as $n \rightarrow \infty$. Since $D$ is closed and $\vec{d}_n \in D$ for all $n$, we have $\vec{d} \in D$. Taking the limit of both sides of 
	$$
	(y,\beta)= \lambda_n \vec{d}_n,
	$$
	as $n \rightarrow \infty$, we have $(y,\beta) = \vec{d}r$ with $\vec{d} \in D$ and $r$ being the attained infimum of all such values such that $(y,\beta) \in rD$. This shows the desired result.
\end{proof}

The following theorem provides conditions that guarantee nonemptiness of the fixed point set. The strategy is to relate an exposed face of the fundamental set $D$ to the fixed points of the algorithm.

\begin{theorem}[Existence of fixed points of \GPPPPA]\label{thm:existence}
	Let $D$ be the (closed) fundamental set of $\kappa$ as in \eqref{eqn:Minkowski}. The following hold.
	\begin{enumerate}[label=(\roman*)]
		\item\label{nonempty:case1} If there exists $\lambda' \geq 0$ such that $$
		\lambda' = \max_{\lambda \in \R}\{\lambda \;\;|\;\; \exists y \in X \text{so\;that}\;(y,\lambda) \in D  \},$$
		then $F=D \cap \{(y,\lambda')\; |\; y \in X\}$ is an exposed face of $D$ and 
		\begin{enumerate}
			\item\label{nonempty:case1a} $\left(\frac{1}{1+\lambda'} \right)F = T(\Fix P_S \circ T);$ and
			\item\label{nonempty:case1b} Any $(x,1) \in \Fix P_S \circ T$ satisfies $T(x,1)=\left(x,\frac{\lambda'}{1+\lambda'}\right)$.
		\end{enumerate}

		For example, this is always the case when $D$ is bounded.
		\item\label{nonempty:case2} If such a $\lambda'$ does not exist and there exists a sequence $(y_n,\lambda_n)_{n \in \N}$ such that $\lambda_n\rightarrow \infty$ and $\lambda_n/\|y_n\| \rightarrow m>0$, then $T(\Fix P_S \circ T) = \Fix P_S \circ T = \zer \kappa \cap S$.
	\end{enumerate}
\end{theorem}
\begin{proof}
	\ref{nonempty:case1}: Suppose $\lambda'$ exists as described. To understand why $F$ is an exposed face of $D$, see Remark~\ref{rem:exposedface} below.
	
	\textbf{Case 1:} $\lambda'=0$, then any point in $F$ is of the form $(x,0)$ for some $x \in X$. Moreover, $(0,1)$ is a fixed point and satisfies $\|T(0,1)-(0,0)\|=1$. It is then a straightforward consequence of Proposition~\ref{prop:fixedpoints} that $T(x,1)=(x,0)$ if and only if $(x,1)$ is a fixed point of $P_S\circ T$. This is all we needed to show in this case.
	
	\textbf{Case 2:} $\lambda'>0$. Since $\lambda'>0$, the set $\cup_{r\geq0}rD\cap (X\times \{\theta \}) \neq \emptyset$ for each $\theta>0$. For any $(x,\theta) \in {\rm cone}D$ with $\theta>0$ we have from Lemma~\ref{lem:conedom}\ref{lem:conedom2} that there exists minimal $r_{(x,\theta)}>0$ such that $(x,\theta) \in r_{(x,\theta)} D.$ By the Minkowski definition of the gauge this means, 
	$$
	\kappa(x,\theta)=r_{(x,\theta)}.
	$$
	There must also exist $(d_x,\mu_\theta) \in D$ such that $r_{(x,\theta)}(d_x,\mu) = (x,\theta)$. Thus $r_{(x,\theta)}\mu_\theta=\theta$ and $r_{(x,\theta)} =\theta/\mu_\theta$. Furthermore, our choice of $\lambda'$ guarantees that $\mu_\theta \leq \lambda'$. Using positive homogeneity we have
	$$
	r_{(x,\theta)} = \kappa(x,\theta) = \kappa(r_{(x,\theta)}(d_x,\mu_\theta)) = r_{(x,\theta)}\kappa(d_x,\mu_\theta)\quad\text{and\;so}\quad \kappa(d_x,\mu_\theta)=1.
	$$
	Again using positive homogeneity, we obtain
	\begin{equation}\label{theta-mu}
	\kappa(x,\theta) = \kappa((\theta/\mu_\theta)(d_x,\mu_\theta))=(\theta/\mu_\theta)\kappa((d_x,\mu_\theta)) =\theta/\mu_\theta,
	\end{equation}
	where the final equality is because we just showed $\kappa(d_x,\mu_\theta)=1$.
	
	Additionally, for any point $(y',\lambda') \in F$, homogeneity assures that
	\begin{equation}\label{theta-lambdaprime}
	\theta/\lambda' \geq (\theta/\lambda') \kappa(y',\lambda')=\kappa\left((\theta/\lambda') (y',\lambda')\right) = \kappa((p_\theta y',\theta)) \quad\text{where}\quad p_\theta:=\theta/\lambda'.
	\end{equation}
	Now let 
	\begin{equation}
	\theta':= \underset{\theta \in \R}{\argmin}\max \{\theta/\lambda',|\theta-1| \}=\frac{\lambda'}{1+\lambda'}.\label{thetaargmin}
	\end{equation}
	Now we show that $(p_{\theta'} y',1)$ is in $\Fix P_S \circ T$. Remember that $\kappa_1$ is the polar envelope of $\kappa$ from Definition~\ref{def:pprox}. We have the following.
	\begin{subequations}\label{theta}
	\begin{align}
	\kappa_1((p_{\theta'} y',1))&=\underset{(x,\theta)\in X\times \R}{\inf}\max\{\kappa(x,\theta),\|(x,\theta)-(p_{\theta'} y',1)\| \} \nonumber\\
	 &= \underset{(x,\theta) \in {\rm cone}D}{\inf}\max\{\kappa(x,\theta),\|(x,\theta)-(p_{\theta'} y',1)\| \} \label{theta2} \\
	&= \underset{(x,\theta) \in {\rm cone}D}{\inf}\max\{\theta/\mu_\theta,\|(x,\theta)-(p_{\theta'} y',1)\| \} \label{theta3} \\
	&\geq \underset{\theta \in \R}{\min}\max\{\theta/\lambda',|\theta-1| \} \label{theta4}\\
	&=\max\{\theta'/\lambda',|\theta'-1| \}\label{theta5} \\
	&\geq \max\{\kappa(p_{\theta'} y',\theta'),\|(p_{\theta'} y',\theta')-(p_{\theta'} y',1)\| \}. \label{theta6}
	\end{align}
	\end{subequations}
	Here \eqref{theta2} is true by Lemma~\ref{lem:conedom}\ref{lem:conedom1}, \eqref{theta3} holds by \eqref{theta-mu}, \eqref{theta4} holds because $\mu_{\theta}\leq \lambda'$, \eqref{theta5} holds by \eqref{thetaargmin}, and \eqref{theta6} is obtained by applying \eqref{theta-lambdaprime} with $\theta = \theta'$. Altogether \eqref{theta} shows that $(p_{\theta'}y',\theta')=T(p_{\theta'}y',1)$, and so $(p_{\theta'}y',1)\in \Fix P_S \circ T$.
	
	Notice that $\theta' \in \left ]0,1\right[$ and is nearer to $1$ for larger $\lambda'$ and nearer to $0$ for smaller $\lambda'$, exactly as we would expect. Notice also that we have shown that \emph{any} point 
	$$
	(\theta'/\lambda')(y',\lambda') = \frac{1}{1+\lambda'}(y',\lambda')\in \frac{1}{1+\lambda'}F,
	$$
	admits a corresponding point $(p_{\theta'}y',1) \in \Fix P_S \circ T$ whose proximal image is
	$$
	T(p_{\theta'}y',1) = (p_{\theta'}y',\theta')=(\theta'/\lambda')(y',\lambda').
	$$
	This shows that 
	$$
	\frac{1}{1+\lambda'}F \subset T\left(\Fix P_S \circ T\right).
	$$
	
	Now let $(x,1) \in \Fix P_S \circ T$. Using Proposition~\ref{prop:fixedpoints}, we have that any $(x,1) \in \Fix P_S \circ T$ must satisfy $T(x,1)=(x,\lambda)$ where $\lambda \in \left[0,1\right]$ and
	$$
	\|(x,1) - (x,\lambda)\| = \|(p_{\theta'} y',\theta')-(p_{\theta'} y',1)\|=r' = |1-\theta'|,
	$$
	which forces $\lambda = \theta'$. This shows that \ref{nonempty:case1b} is true. 
		
	Now by Lemma~\ref{lem:conedom}\ref{lem:conedom1},$(x,\theta') \in \cone D$, since $(x,\theta') \in \dom \kappa$. Now using the fact that $(x,\theta')\in \cone D$ and following the same reasoning as we used above to obtain \eqref{theta-mu}, we have that there exists a value $r_{(x,\theta')}$ and a point $(d_x,\mu_{\theta'})$ with $\mu_{\theta'} \leq \lambda'$ such that $r_{(x,{\theta'})}(d_x,\mu_{\theta'}) = (x,{\theta'})$ and $\kappa(x,{\theta'})=({\theta'}/\mu_{\theta'})$. Since $(x,1),(p_{\theta'}y',1) \in \Fix P_S \circ T$, we have from Proposition~\ref{prop:fixedpoints} that 
	$$
	|{\theta'}-1|=\|(x,{\theta'})-(x,1)  \|= \|(p_{\theta'} y',\theta')-(p_{\theta'} y',1)\|=r'.
	$$
	Moreover, ${\theta'}/\mu_{\theta'}=\kappa(x,{\theta'}) \leq \|(x,{\theta'})-(x,1)\|=r'$, and so
	\begin{align}
	r'&=\max \{{\theta'}/\mu_{\theta'},|1-{\theta'}|\}, \nonumber\\
	&\geq \underset{(\theta,\mu) \in \R_+ \times \left]0,\lambda'\right]}{\min}\max \{\theta/\mu,|1-\theta|\}, \nonumber\\
	&=\max \{\theta'/\lambda',|1-\theta'|\},\nonumber \\
	&=r'. \label{lambdamustbetheta}
	\end{align}
	The equality throughout \eqref{lambdamustbetheta} forces $\mu_{\theta'} = \lambda'$. Finally,	
	\begin{align*}
	(d_x,\mu_{\theta'})&= (d_x,\lambda') \in F,\\
	\text{and\; so}\quad T(x,1)&= (x,{\theta'}) = (x,\theta') = (\theta'/\lambda')(d_x,\lambda') \in (\theta'/\lambda')F = \frac{1}{1+\lambda'}F.
	\end{align*}
	This shows that
	$$
	\frac{1}{1+\lambda'}F \supset T\left(\Fix P_S \circ T\right).
	$$
	This concludes the proof of \ref{nonempty:case1a}.
	
	\ref{nonempty:case2}: Let the sequence $(y_n,\lambda_n)_{n \in \N}$ exist as described. By compactness of the unit ball in Euclidean space and by appealing to a subsequence if necessary, the sequence $y_n/\|y_n\|$ converges to some $y$ in the unit ball in $X$. Now since $(y_n,\lambda_n) \in D$ for all $n$,
	$$
	\kappa(y_n,\lambda_n) \leq 1\quad (\forall n),
	$$
	and by the Minkowski function representation of $\kappa$,
	$$
	\kappa\left(\frac{y_n}{m\|y_n\|},\frac{\lambda_n}{m\|y_n\|}\right) = \frac{1}{m\|y_n\|}\kappa\left(y_n,\lambda_n\right).
	$$
	Taking the limits of both sides as $n\rightarrow \infty$ and using the lower semicontinuity of $\kappa$, we obtain
	$$
	\kappa\left(\frac{y}{m},1\right) = 0.
	$$
	The point $\left(\frac{y}{m},1\right) \in \zer \kappa \cap S$ is clearly a fixed point of $T$ since
	$$
	\max \left\{\kappa\left(\frac{y}{m},1\right),\left\|\left(\frac{y}{m},1\right)-\left(\frac{y}{m},1\right)\right\| \right\} = 0.
	$$
	Thereafter appealing to Lemma~\ref{lem:FixT}, Proposition~\ref{prop:fixedpoints}, and the fact that $\kappa(T(x,1)) \leq \|(x,1)-T(x,1)\|=r'=0$ for any $(x,1) \in \Fix P_S \circ T$, the result \ref{nonempty:case2} is clear.
\end{proof}

\begin{remark}[What do we mean by an \emph{exposed face}?]\label{rem:exposedface}
	Let us explain what we mean in Theorem~\ref{thm:existence} when we say that $F$ is an exposed face of $D$. Recalling \cite[Definition~6]{roshchinafaces}, $F$ is an exposed face of a closed, convex set $D$ if there exists a supporting hyperplane $H$ to $D$ with $F = D \cap H$. In our case, $H = X \times \{\lambda' \}$. Recalling \cite[Definition~5]{roshchinafaces}, the hyperplane $H$ is a \emph{supporting hyperplane} because $D$ lies entirely in the affine half space $X \times \R_{\leq \lambda'}$ defined by $H$.
\end{remark}

The following example showcases a situation when $\Fix P_S \circ T$ may be empty. In so-doing, it illustrates the importance of the condition $m>0$ in Theorem~\ref{thm:existence}\ref{nonempty:case2}.

\begin{example}\label{ex:whyismgreaterthanzero}
	Let $D= \{(y,\lambda)\; |\; y \geq \lambda^2 \} \subset \R^2$. Then for any $(y,1) \in S$, $T(y,1)=(u,1)$ with $u>y$, and so $\Fix P_S \circ T = \emptyset$.
\end{example}

\subsection{Fixed points of \PPPPA: facial characterization}\label{sec:facial}

When we take the results of Theorem~\ref{thm:existence} and specify from $\kappa$ back to the perspective transform $f^\pi$, we recover the following characterization of the fixed points of \PPPPA.

\begin{theorem}[Facial characterization of fixed points of \PPPPA ]\label{thm:fixedpoints_argmin_new}
	Let $f:X \rightarrow \R_+ \cup \{+\infty \}$ be a proper closed nonnegative convex function with $\inf f>0$ and $\argmin f \neq \emptyset$. Let $\kappa = f^\pi$. The following hold:
	\begin{enumerate}[label=(\roman*)]
		\item\label{fixed1} $\lambda' = \frac{1}{\min_{u \in X}f(u)},$ where $\lambda'$ is as in Theorem~\ref{thm:existence}\ref{nonempty:case1};
		\item\label{fixed2} Where $F$ is as in Theorem~\ref{thm:existence}\ref{nonempty:case1}, $F=\frac{1}{\min_{u \in X}f(u)}(\argmin f \times \{1\})$;
		\item\label{fixed3} Any $(x,1) \in \Fix P_S \circ T$ satisfies $\frac{1+\lambda'}{\lambda'}x \in \argmin f$;
		\item\label{fixed4} Any $y \in \argmin f$ satisfies $(\frac{\lambda'}{1+\lambda'}y,1) \in \Fix P_S \circ T$.
	\end{enumerate}
\end{theorem}
\begin{proof}
	\ref{fixed1}: For simplicity, let $\eta:=\min_{u \in X}f(u)$. We will first show that 
	$$
	\frac{1}{\eta} = \max_{\lambda \in \R} \{\lambda \;\;|\;\; \exists y \in X\; \text{so\;that}\;(y,\lambda) \in D \}.
	$$
	Let $y \in \argmin f$. Then
	$$
	f^\pi(y/\eta,1/\eta) = \frac{1}{\eta} f^\pi(y,1) =\frac{1}{\eta}f(y)=\eta = 1,
	$$
	and so $(y,1/\eta) \in D$. To see that $1/\eta$ is maximal, suppose for a contradiction that there exists $(y_0,\lambda_0) \in \left( X \times \left]1/\eta,\infty \right[\right) \cap D.$ Then 
	$$
	1 \geq f^\pi(y_0,\lambda_0) \geq  \lambda_0 f^\pi(y_0/\lambda_0,1) > \frac{1}{\eta}f(y_0/\lambda_0,1) = \frac{1}{\eta}f(y_0/\lambda_0) \geq 1,
	$$
	a contradiction.
	
	\ref{fixed2}: Having shown \ref{fixed1}, we have from the definition of $F$ that $F=D \cap \left\{(x,1/\eta)\;|\; x \in X \right\}$. Let $(y,1/\eta) \in F$. Then
	$$
	1 \geq f^\pi(y,1/\eta)=(1/\eta) f^\pi(y\eta,1) = (1/\eta) f(y\eta) \geq 1,
	$$
	and the equality throughout forces $f(y\eta)=\eta$. Thus $y\eta \in \argmin f$ and so $y \in (1/\eta)\argmin f$. Thus $F \subset (1/\eta)\argmin f$. The reverse inclusion is similar.
	
	\ref{fixed3} \& \ref{fixed4}: By Theorem~\ref{thm:existence}\ref{nonempty:case1a} $(x,1) \in \Fix P_S \circ T$ is equivalent to 
	\begin{equation}\label{eqn:fixed34}
	\left(x, \frac{\lambda'}{1+\lambda'} \right) = \left(\frac{1}{1+\lambda'} \right) \left(y,\frac{1}{\eta}\right)\;\;\text{for\;some}\;\;(y,1/\eta) \in F.
	\end{equation}
	Having shown \ref{fixed2}, we have that the latter inclusion is equivalent to $y \eta \in \argmin f$. Combining with \eqref{eqn:fixed34}, 
	$$
	y\eta \in \argmin f \iff \left( \frac{1+\lambda'}{1}x\right) \eta \in \argmin f.
	$$
	Having shown \ref{fixed1}, this is equivalent to
	$$
	\frac{1+\lambda'}{\lambda'}x \in \argmin f,
	$$
	which shows both \ref{fixed3} and \ref{fixed4}.
	
	%\ref{fixed3}: Applying Theorem~\ref{thm:fixedpoints_argmin}\ref{fixedpoints1}, we have that $(x,\lambda_*) = T(x,1)$ satisfies $\lambda_*^{-1}x \in \argmin f$. Additionally, because we just showed \ref{fixed1}, we know from Theorem~\ref{thm:existence}\ref{nonempty:case1} that $\lambda_* = \frac{\lambda'}{1+\lambda'}$, which shows the result. 
	
	%\ref{fixed4}: From Theorem~\ref{thm:fixedpoints_argmin}\ref{fixedpoints2}, $((1+\eta)^{-1} y,1) \in \Fix P_S \circ T$. From \ref{fixed1}, $(1+\eta)^{-1} = \frac{\lambda'}{1+\lambda'}$, which shows the result.
\end{proof}

In the following remark, we compare the facial characterization of fixed points of \PPPPA \; from Theorem~\ref{thm:fixedpoints_argmin_new} with the closely related results from \cite{friedlander2019polar}.

\begin{remark}[On synchronicity between Theorems~\ref{thm:fixedpoints_argmin} and \ref{thm:fixedpoints_argmin_new}]\label{rem:synchronicity}
	Theorem~\ref{thm:fixedpoints_argmin_new} subsumes and is closely connected with the original results of {\cite[Theorem~7.4]{friedlander2019polar}}, which we recalled as Theorem~\ref{thm:fixedpoints_argmin}. To see why, notice that items \ref{fixed1}, \ref{fixed3}, \ref{fixed4} of Theorem~\ref{thm:fixedpoints_argmin_new} have the following characterizations.
	
	\ref{fixed3}: Applying Theorem~\ref{thm:fixedpoints_argmin}\ref{fixedpoints1}, we have that $(x,\lambda_*) = T(x,1)$ satisfies $\lambda_*^{-1}x \in \argmin f$. Theorem~\ref{thm:existence}\ref{nonempty:case1} guarantees that $\lambda_* = \frac{\lambda'}{1+\lambda'}$.
	
	\ref{fixed4}: From Theorem~\ref{thm:fixedpoints_argmin}\ref{fixedpoints2}, $((1+\eta)^{-1} y,1) \in \Fix P_S \circ T$. From Theorem~\ref{thm:fixedpoints_argmin_new}\ref{fixed1}, $(1+\eta)^{-1} = \frac{\lambda'}{1+\lambda'}$.
	
	\ref{fixed1}: The condition $1/(1+\eta) = \frac{\lambda'}{1+\lambda'}$ then yields $\lambda'=1/\eta$.
	
	Theorem~\ref{thm:fixedpoints_argmin_new} essentially uses the more general results from Theorem~\ref{thm:existence} to show that the minimizers of $f$ form an exposed face of $(\min_{u \in X}f(u))D$: namely the face that is $(\min_{u \in X}f(u))F$. The other items are all a natural consequence of this.
\end{remark}

\section{Conclusion}\label{sec:conclusion}

We now state our eponymous convergence result, which shows global convergence of \PPPPA \;in the full generality of \cite{friedlander2019polar}. It also, under sufficient conditions to guarantee existence of a fixed point, shows convergence of \GPPPPA.

\begin{theorem}[Convergence of \PPPPA \; and \GPPPPA]\label{thm:convergence}
	Let $D$ be the (closed) fundamental set of $\kappa$ as in \eqref{eqn:Minkowski}. Suppose one of the following holds.
	\begin{enumerate}[label=(\roman*)]
		\item\label{convergence:case4} $\kappa = f^\pi$ for $f:X \rightarrow \R_+ \cup \left\{+\infty \right\}$ a proper closed nonnegative convex function with $\inf f=0$ and $\argmin f \neq \emptyset$; 
		\item\label{convergence:case1} $\kappa = f^\pi$ for $f:X \rightarrow \R_+ \cup \left\{+\infty \right\}$ a proper closed nonnegative convex function with $\inf f>0$ and $\argmin f \neq \emptyset$; 
		\item\label{convergence:case2} There exists $\lambda' \geq 0$ such that $$
		\lambda' = \max_{\lambda \in \R}\{\lambda \;\;|\;\; \exists y \in X \text{so\;that}\;(y,\lambda) \in D  \};$$
		\item\label{convergence:case3} Such a $\lambda'$ does not exist and there exists a sequence $(y_n,\lambda_n)_{n \in \N}$ such that $\lambda_n\rightarrow \infty$ and $\lambda_n/\|y_n\| \rightarrow m>0$.
	\end{enumerate}
	Let $\gamma \in \left[0,1 \right[$ and $(y_0,1) \in S$. Then the following hold.
	\begin{enumerate}
	\item The sequence given by 
	\begin{align*}
	(y_{n+1},1):=&\;\mathcal{U}_\gamma (y_{n},1), \nonumber\\
	\text{where}\quad \mathcal{U}_\gamma:=&\;(1-\gamma)P_S \circ T +\gamma \Id %\label{relaxedorunrelaxed}
	\end{align*}
	is convergent to some $(y,1) \in \Fix P_S \circ T$;
	
	\item The shadow sequences $(y_{n+1},\lambda_{n+1})= Tx_n$ satisfy $\lambda_n \rightarrow \lambda$ for some $\lambda \in \left[0,1\right]$; 
	
	\item When \ref{convergence:case1} or \ref{convergence:case2} holds, $\lambda= \frac{\lambda'}{1+\lambda'}$; 
	
	\item When \ref{convergence:case1} holds, $\lambda' = 1/(\inf f)$ and $\left(\frac{1}{\lambda_n}\right) y_n \rightarrow \left(\frac{1+\lambda'}{\lambda'}\right) y \in \argmin f$;
	
	\item When \ref{convergence:case4} holds, $y_n \rightarrow y \in \argmin f$ and $\lambda_n \rightarrow 0$.
	\end{enumerate}
\end{theorem}
\begin{proof}
	\; \\
	\textbf{Fixed points:} \ref{convergence:case4}: Since $\inf f=0$, any $x \in \argmin f$ satisfies
	$$
	(f^\pi(x,1)=0) \underset{(Lemma~\ref{lem:FixT})}{\implies} ((x,1) \in \Fix T) \implies ((x,1) \in \Fix P_S \circ T).
	$$
	By Theorem~\ref{thm:fixedpoints_argmin_new}, we have that \ref{convergence:case1} $\implies$ \ref{convergence:case2}. Either of the assumptions \ref{convergence:case2} or \ref{convergence:case3} guarantees existence of a fixed point of $P_S \circ T$ by Theorem~\ref{thm:existence}. 
	
	\noindent \textbf{Convergence:} Having shown that a fixed point exists, convergence of $(y_n)_n$ is assured by Corollary~\ref{cor:relaxedandunrelaxed}, and the convergence of $(\lambda_n)_n$ is guaranteed by Corollary~\ref{cor:lambdashadow}. The characterization of $\lambda'$ in cases \ref{convergence:case1} and \ref{convergence:case2} is due to Theorems~\ref{thm:existence} and \ref{thm:fixedpoints_argmin_new}.
\end{proof}

\subsubsection*{Further research}

We suggest three further avenues of inquiry. Firstly, results on faces of fundamental sets (e.g. Theorem~\ref{thm:fixedpoints_argmin_new}) are of interest in the development of more general theory. Secondly, Friedlander, Mac\^{e}do, and Pong also introduced a second algorithm, $\mathbf{EMA}$, which is not addressed here \cite{friedlander2019polar}. A natural question is whether $\mathbf{EMA}$ possesses similar properties to \PPPPA . Finally, a motivating question is whether or not algorithms such as \PPPPA\; may have computational advantages for certain problems.

\subsubsection*{Acknowledgements}

The author was supported by Hong Kong Research Grants Council PolyU153085/16p. The author thanks Ting Kei Pong and Michael P. Friedlander for their useful suggestions on this manuscript. % and by the Alf van der Poorten Traveling Fellowship (Australian Mathematical Society).

\subsubsection*{Data Availability Statement}

Data availability considerations are not applicable to this research.

\bibliographystyle{plain}
\bibliography{bibliography}

\end{document}